\documentclass[11pt]{amsart}

\usepackage[all]{xy}
\def\pf{\begin{proof}}
\def\epf{\end{proof}}

\usepackage[active]{srcltx}   

\newcommand{\supp}{\operatorname{sup}}
\newcommand{\Aut}{\operatorname{Aut}}
\newcommand{\Hom}{\operatorname{Hom}}
\newcommand{\End}{\operatorname{End}}

\newcommand{\Endd}{\mathcal{END}}

\newcommand{\Stab}{\operatorname{Stab}}
\newcommand{\Vect}{\operatorname{Vect}}
\newcommand{\vect}{\operatorname{vect}}
\newcommand{\id}{\operatorname{id}}
\newcommand{\co}{\operatorname{co}}

\newcommand{\nd}{\noindent}
\newcommand{\ot}{\otimes}
\newcommand{\vu}{\vspace{.1cm}}

\providecommand{\bysame}{\makebox[3em]{\hrulefill}\thinspace}
\newcommand{\F}{\Bbbk}

\newcommand{\E}{\mathcal{E}}
\newcommand{\Hc}{\mathcal{H}}
\newcommand{\D}{\mathcal{YD}}
\newcommand{\I}{{\mathbb I}}

\newcommand{\q}{\mathbf q}
\newcommand{\Z}{\mathbb Z}
\newcommand{\N}{\mathbb N}

\newcommand{\Go}{\Gamma_0}
\newcommand{\Lo}{\Lambda_0}

\newcommand{\trid}{\triangleright}
\newcommand{\trii}{\triangleleft}
\newcommand{\ydg}{^{\F\Gamma}_{\F\Gamma}\mathcal{YD}}
\newcommand{\ydl}{^{\F\Lambda}_{\F\Lambda}\mathcal{YD}}
\newcommand{\ydG}{^{\F G}_{\F G}\mathcal{YD}}
\newcommand{\zt}{\mathbb{Z}^{\theta}}
\newcommand{\A}{\mathcal{A}}
\newcommand{\B}{\mathcal{B}}
\newcommand{\C}{\mathcal{C}}
\newcommand{\Dc}{\mathcal{D}}
\newcommand{\Vc}{\mathcal{V}}
\newcommand{\cop}{\operatorname{cop}}

\usepackage{amsmath,amsthm,amssymb,amsxtra,amscd, amsbsy,eucal,color}

\newtheorem{theorem}{Theorem}[section]
\newtheorem{lema}[theorem]{Lemma}

\newtheorem{prop}[theorem]{Proposition}

\theoremstyle{definition}
\newtheorem{definition}[theorem]{Definition}
\newtheorem{exa}[theorem]{Example}

\theoremstyle{remark}
\newtheorem{remark}[theorem]{Remark}

\numberwithin{equation}{section}

\begin{document}

\title[Examples of pointed color Hopf algebras]{Examples of pointed color Hopf algebras}
\author[Andruskiewitsch, Angiono,  Bagio]
{Nicol\'as Andruskiewitsch, Iv\'an  Angiono, Dirceu Bagio}

\address{FaMAF-Universidad Nacional de C\'ordoba, CIEM (CONICET),
Medina A\-llen\-de s/n, Ciudad Universitaria, 5000 C\' ordoba, Rep\'ublica Argentina.} \email{(andrus|angiono)@famaf.unc.edu.ar}

\address{Departamento de Matem\'atica, Universidade Federal de Santa Maria,
97105-900, Santa Maria, RS, Brazil} \email{bagio@smail.ufsm.br}

\thanks{\noindent 2000 \emph{Mathematics Subject Classification.}
16T05,17B75. \newline
N. A. and I. A. were partially supported by CONICET, FONCyT-ANPCyT, Secyt (UNC).
D. B. was supported through CNPq - 201586/2010-0, CAPES-MATH-AmSud}

\begin{abstract}
We present examples of color Hopf algebras, i.~e. Hopf algebras in  color categories (braided tensor categories with braiding induced
by a bicharacter on an abelian group), related with quantum doubles of pointed Hopf algebras.
We also discuss semisimple color Hopf algebras.
\end{abstract}

\maketitle

\section*{Introduction}

Color Lie algebras appeared in \cite{Ree} (under the name \emph{generalized Lie algebra of type $\chi$}), were rediscovered in
\cite{Sc0} and then studied in several papers, e.~g. \cite{BD,Hu,Kory, chen1, Sc2,chen, Mo}.
A color Lie algebra is a Lie algebra in a symmetric tensor category, namely the category of $G$-graded vector spaces over a finite abelian group $G$ with symmetry
given by a skew-symmetric bicharacter $\beta:G\times G\to \F^\times$.
Hopf algebras, unlike Lie algebras, can be defined in braided tensor categories. In this paper we study
Hopf algebras in the braided tensor category of $G$-graded vector spaces with braiding
given by a (not necessarily skew-symmetric) bicharacter $\beta:G\times G\to \F^\times$: we call them \emph{color Hopf algebras}.
Pointed Hopf superalgebras were considered in \cite{AAY} in relation with the classification of
finite-dimensional Nichols algebras of diagonal type \cite{He}. However
some of the Nichols algebras in the list of \cite{He} are neither of standard (close to simple Lie algebras) nor of super type.
We wonder whether they are related with color Hopf algebras, and this is one of the motivations of the present article,
whose contents we describe next.
Section \ref{section:preliminaries} is devoted to basic constructions and results in color categories; albeit
these are particular instances of constructions and results in braided tensor categories, we feel that the explicit
formulae in this context might be useful for the reader.
We also relate the symmetric color categories with the classification of semisimple triangular Hopf algebras over a algebraically closed field
\cite{EtGe, EtGe1}.
In Section \ref{section:color categories}, we adapt the lifting method \cite{AS3} to the setting of color Hopf algebras.
In Section \ref{section:drinfeld doubles}, we produce color Hopf algebras versions of
the quantum doubles of bosonizations of color Nichols algebras.
We introduce the notion of \emph{consistent} coloring of a generalized Dynkin diagram,
meaning that is preserved by the movements of the Weyl groupoid.
In Section \ref{section:semisimple color}, we discuss some examples of semisimple color Hopf algebras.

\section{Preliminaries and notation}\label{section:preliminaries}

Let $\F$ be an algebraically closed field of characteristic 0 and $\F^{\times}=\F-\{0\}$;
all vector spaces, algebras and tensor products are over $\Bbbk$.

Throughout this paper, $G$ is a finite abelian group, denoted multiplicatively, and $A = \widehat{G} = \Hom(G,\F^{\times})$
is the character group of $G$. Then $(e_g)_{g\in G}$ denotes the canonical base of $\F G$ and $\Aut G$ the group of automorphisms of $G$.
The elements of the dual basis of $(\F G)^* =\F^G$ are denoted by $\delta_g$, $g\in G$.
Let $L$ be an abelian group and suppose that $G$ acts on $L$; we denote by $Z^1(G,L)$ (resp. $Z^2(G,L)$) the set of normalized
1-cocycles (resp. 2-cocycles).
If $A$ acts on a set $X$ and $x\in X$, then $\Stab_A(X)$ is the stabilizer of the action and $A^x$ is the isotropy
subgroup of $x$. $C_N$, $N\in\N$, is the cyclic group of order $N$. For an arbitrary group $\Gamma$, $Z(\Gamma)$ denotes its center. 

We also denote by $G(C)$ the set of all group-like elements of a coalgebra $C$. The category of left, respectively right, $C$-comodules is denoted by $^{C}\mathcal{M}$; resp. $\mathcal{M}^{C}$. If $B$ is an algebra, then the category of left, respectively right, $B$-modules is denoted by $_{B}\mathcal{M}$;
resp. $\mathcal{M}_{B}$. When $B$ is a bialgebra, $\mathcal{P}(B)$ denotes the set of all primitive elements.

All Hopf algebras in this paper have bijective antipodes.
Let $H$ be a Hopf algebra. We will denote by $\Aut H$ the group of Hopf algebra automorphisms
of $H$; and by $\mathcal{HZ}(H)$ the Hopf center of $H$, i.~e., the maximal central Hopf subalgebra of $H$ \cite{A}.
We assume that the reader has familiarity with tensor categories, see e.~g. \cite{BK}, and
with Yetter-Drinfeld modules and bosonization, see e.~g. \cite{AS3}. Recall that a tensor functor between  tensor categories $\C$ and $\Dc$ (with associativity $a$) is a pair $(F, \phi)$,
where $F:\C \to \Dc$ is an additive $\F$-linear functor and
$\phi_{X,Y}: F(X\ot Y) \to F(X)\ot F(Y)$ is a natural isomorphism such that
$a_{F(X), F(Y), F(Z)} \left(\phi_{X,Y} \ot \id\right)\phi_{X\ot Y, Z} = \left(\id \ot \phi_{Y,Z} \right)
\phi_{X, Y\ot Z} F(a_{X, Y, Z})$ for all $X,Y,Z\in \C$, and some compatibility for the unit.
If $\C$ and $\Dc$ are braided
(with braiding $c$), a monoidal tensor functor $F:\C \to \Dc$ is braided when for all $X,Y\in \C$
\begin{align}\label{eq:mon-tensor-functor}
c_{F(X), F(Y)}\phi_{X,Y} &= \phi_{Y,X} F(c_{X,Y}).
\end{align}

\subsection{Color categories}

We denote by $\Vect^{G}$ the category of $G$-graded vector spaces.
Let $V=\oplus_{g\in G} V_g\in \Vect^G$ and $v\in V$. If $v\in V_g$, then we write $|v|=g$.
The support of $V$ is $\text{sup}\,V=\{g\in G\,:\,V_g\neq 0\}$. If $v=\sum_{g\in G}v_g$, then
the support of $v$ is $\text{sup}\,v=s(v)=\{g\in G\,:\,v_g\neq 0\}$.
The category $\Vect^G$ is a tensor one: if $V,W\in \Vect^{G}$, then $V\otimes W\in \Vect^{G}$
with $(V\otimes W)_g=\oplus_{hl=g}V_h\otimes W_l$.
Let $\F_g\in \Vect^G$, $g\in G$ be given by $(\F_g)_l=\F$ if $l=g$ and $(\F_g)_l=0$ otherwise.
Thus any object is a direct sum of copies of various $\F_g$.
Then $\F_g\otimes\F_h\simeq \F_{gh}$. The unit object is  $\F=\F_1$.
The full subcategory $\vect^G$ of finite dimensional $G$-graded vector spaces is rigid: if $V\in \vect^G$,
then $V^{*}=\Hom(V,\F)=\oplus_{g\in G}V^{*}_g$, with $V^{*}_g=\Hom(V_{g^{-1}},\F)$.

\medbreak
A {\it bicharacter} on $G$ is a $\mathbb{Z}$-bilinear form  $\beta:G\times G\to \Bbbk^{\times}$, that is,
\begin{align}\label{color:bilineal}
&\beta(gh,l)=\beta(g,l)\beta(h,l),& &\beta(g,hl)=\beta(g,h)\beta(g,l),& &g,h,l\in G.&
\end{align}
A {\it commutation factor} on $G$ is a skew-symmetric bicharacter $\beta$, that is

\begin{align}\label{color:commutation}
&\beta(g,h)\beta(h,g)=1,& g,h\in G.&
\end{align}
Let  $\beta$ be a bicharacter on $G$.
Then $\Vect^{G}$ is braided with braiding given by
\begin{align*}
c & = c_{\beta}:V\otimes W\to W\otimes V, & c(v\otimes w) &= \beta(g,h)w\otimes v, & v&\in V_g,\quad w\in W_h.
\end{align*}
Here it is crucial that $G$ is abelian, for $c$ to be a map of $G$-graded vector spaces.
Conversely, any braiding in the tensor category  $\Vect^{G}$ comes from a  bicharacter $\beta$ on $G$.
Furthermore, $c_{\beta}$ is a symmetry in $\Vect^{G}$ if and only if  $\beta$ is a commutation factor on $G$.

\begin{definition}\label{def:color category} Let $G$ be a finite abelian group and $\beta$ a bicharacter on $G$.
We denote by $\Vect^{G}_\beta$ the braided tensor category $\Vect^{G}$ with braiding $c_{\beta}$.
Any braided tensor category like this  is called a color category with color $\beta$.
\end{definition}

\begin{remark}\label{rem:color-categories}
(i). Let $G$, $\Gamma$ be finite abelian groups,  $\beta$ and $\vartheta$ bicharacters on $G$ and $\Gamma$
respectively, $\psi: G \to\Gamma$ a group homomorphism and $\gamma:G\times G\to \F^{\times}$ a normalized 2-cocycle.
Then there is  a tensor functor $(F, \phi): \Vect^{G} \to \Vect^{\Gamma}$ given by
$F(\F_g) = \F_{\psi(g)}$, $\phi_{\F_g, \F_h} = \gamma(g,h)$.
Assume that for all $g,h\in G$
\begin{align}\label{eq:mon-tensor-functor-cociclo}
\beta(g,h) &= \vartheta (\psi(g), \psi(h))\dfrac{\gamma(g,h)}{\gamma(h,g)}.
\end{align}
Then $(F, \phi): \Vect^{G}_\beta \to \Vect^{\Gamma}_\vartheta$ is a braided tensor functor, see \eqref{eq:mon-tensor-functor}.

\medbreak (ii) Let $\vartheta: (A\times G) \times (A\times G) \to \F^{\times}$ be given by
$\vartheta((\xi, g), (\eta, h)) = \eta(g)$. Then  $\Vect^{A\times G}_{\vartheta}$
and $^{\F G}_{\F G}\D$ are equivalent braided tensor categories.

\medbreak (iii)
The color category $\Vect^G_\beta$ is actually isomorphic to a braided subcategory of
$^{\F G}_{\F G}\D$. Indeed, let $\chi, \chi^o: G \to A$ be given by
\begin{align}\label{equation0}
\chi_g(h)=\beta(h,g) & = \chi^o_h(g),&  g,h &\in G.
\end{align}

Then the morphism $G\to A\times G$, $g\mapsto(\chi_g,g)$, together the trivial cocycle, induce a braided
tensor functor $\Vect^G_\beta\to \Vect^{A\times G}_\vartheta\simeq {}^{\F G}_{\F G}\D$.
Explicitly, if $V\in\Vect^{G}$,
then $V\in {}^{\F G}_{\F G}\D$ with action and coaction
\begin{align}\label{eq:explicit-tensorfunctor}
e_g\cdot v  &=\beta(g,|v|)\,v, & \delta(v) &= e_{|v|}\otimes v,&  v&\in V,\, g\in G.
\end{align}
\end{remark}

In what follows, we fix a bicharacter $\beta$ and use the adjective \emph{graded} for notions in $\Vect^G$,
that do not require the braiding,
and  \emph{color} for  notions in $\Vect^G_\beta$, that do.

\subsection{Graded algebras}

A $G$-graded algebra is the same as an associative algebra in $\Vect^{G}$. If $V \in\vect^{G}$, then
$\End (V)$ is a graded algebra with $\End (V)_g=\{T\in \End (V)\,:\,T(V_h)\subset V_{gh},
\,\text{for all}\,\,h\in G\},$
for all $g\in G$.
A {\it graded representation} of a $G$-graded algebra $B$ on $V\in \vect^G$ is a
$G$-graded morphism $\rho:B\to \End (V)$; for general $V$, a graded representation
is a morphism into the subalgebra $\Endd (V) = \oplus_{g\in G}\End (V)_g$
of $\End(V)$.  Note that $\rho:B\to \End (V)$ is a graded
representation if and only if there is an action $B\otimes V\to V$ that respects the grading.
In this case, we say that $V$ is a {\it left graded module}. Analogously, there are right graded modules.
A {\it graded bimodule} is a bimodule such that both the left and the right actions are homogeneous.
The respective abelian categories are denoted by $_B\mathcal{GM}$, $\mathcal{GM}_B$ and
$_B\mathcal{GM}_B$, with obvious morphisms.

\begin{remark}\label{rem:1.1} Let $B$ be a graded algebra and $B_G=B\#\F G$ the smash product algebra;
this is the vector space  $B\otimes \F G$ with the multiplication given by
\begin{align}\label{1.1}
    (x\# e_g)(y\# e_h)& = \beta(g,|y|)xy\# e_{gh},& x,y&\in B,& g,h&\in G.
\end{align}
 Note that if $V\in\, _B\mathcal{GM}$ then $V\in\, _{B_G}\mathcal{M}$ by
 \begin{align}\label{1.2}
    (x\#e_g)v& = \beta(g,|v|)xv, &  x&\in B,& g&\in G, &  v&\in V.
\end{align}
The reciprocal holds provided that $\beta$ is non degenerate. Indeed, if $\xi\in A$, then there exists a unique $h_{\xi}\in G$ such that $\xi(g)=\beta(g,h_{\xi})$, for all $g\in G$. Let $V\in\, _{B_G}\mathcal{M}$.
Since $G$ acts in $V$, $V=\oplus_{\xi\in A}V_{\xi}$. Hence, $V$ is $G$-graded with $V_{h_{\xi}}=V_{\xi}$.
The action of $B$ is the restriction of the action of $B_G$, so $V\in\,_B\mathcal{GM} $.
In this way, the abelian categories $_B\mathcal{GM}$ and $_{B_G}\mathcal{M}$ are equivalent.
\end{remark}

\subsection{Graded coalgebras}

A $G$-graded coalgebra is a coalgebra $C$ in  $\Vect^G$. For instance, the dual coalgebra of a finite dimensional
graded algebra is a graded coalgebra.
A {\it (left) graded comodule} over a graded coalgebra $C$ is a left comodule $V$ such that $V\in \Vect^G$ and
the coaction $\lambda:V\to C\otimes V$ is homogeneous. Right graded comodules and graded bicomodules are defined
in a similar way. These categories are denoted by $^C\mathcal{GM}$, $\mathcal{GM}^C $ and $^C\mathcal{GM}^C$.
Given a graded coalgebra $C$, we denote the coproduct by the following variant of Sweedler's notation:
$\Delta(c)=c^{(1)}\otimes c^{(2)}$,  $c\in C$. Analogously, for a graded comodule $V$,
we write: $\lambda(v)=v^{(-1)}\otimes v^{(0)}$, for all $v\in V$.

\begin{remark}\label{rem:1.2} Let $C$ be a graded coalgebra and $C_G=C\#\F G$ the smash product
coalgebra, i.~e. the vector space  $C\otimes \F G$ with  coproduct and counit
\begin{align}\label{1.3}
    \Delta(c\# e_g)&=c^{(1)}\#e_{|c^{(2)}|g}\otimes c^{(2)}\#e_g, & \varepsilon(c\#e_g) &=\varepsilon(c),& c&\in C, & g&\in G.
\end{align}
As a consequence, $C_G$ is pointed if and only if
 $C$ is pointed and the grading of the coradical $C_0$ is trivial.
In this case,  we say that $C$ is \emph{strictly pointed}.

The abelian categories $^{C}\mathcal{GM}$ and $^{C_G}\mathcal{M}$ are equivalent.
In fact, if $V\in\, ^C\mathcal{GM}$, then $V\in\, ^{C_G}\mathcal{M}$ by
\begin{align}\label{1.4}
    \lambda_G(v) &= v^{(-1)}\#e_{|v^{(0)}|}\otimes v^{(0)}, & v&\in V.
\end{align}
\end{remark}

\subsection{Color Hopf algebras}

The bicharacter $\beta$ is not needed to define algebras and coalgebras in $\Vect^G$, but to define bialgebras
and Hopf algebras. First, $\beta$ allows us to define the twisted product between $G$-graded algebras:
if $B,B'$ are $G$-graded algebras, then $B\underline{\otimes} B'$ is the  $G$-graded algebra
with multiplication $(x\otimes x')(y\otimes y')=\beta(|x'|,|y|)xy\otimes x'y'$, for all $x,y\in B$, $x',y'\in B'$.
A {\it color bialgebra} is a bialgebra $B$ in the category $\Vect^{G}$, that is, $B=\oplus_{g\in G}B_g$ is a $G$-graded
algebra and coalgebra  such that $\Delta:B\to B\underline{\otimes} B$ and $\varepsilon$ are morphisms of $G$-graded algebras.
A color bialgebra $H$ is said a {\it color Hopf algebra} when the identity map has a convolution
inverse $\mathcal{S}\in\End (H)$ (that we shall always assume is bijective.); $\mathcal{S}$ is called the {\it antipode} and satisfies

\begin{enumerate}\renewcommand{\theenumi}{\roman{enumi}}   \renewcommand{\labelenumi}{(\theenumi)}
\item $\mathcal{S}(xy)=\beta(|x|,|y|) \, \mathcal{S}(y)\mathcal{S}(x)$ and

\item $\Delta(\mathcal{S}(x))=\beta(|x^{(1)}|,|x^{(2)}|) \, \mathcal{S}(x^{(2)})\otimes\mathcal{S}(x^{(1)})$,
for all $x,y\in H$.
\end{enumerate}

Evidently, if $\beta = \varepsilon$ is the trivial bicharacter, then a $(G, \varepsilon)$-color Hopf algebra is a $G$-graded Hopf algebra, that is a
Hopf algebra which is $G$-graded algebra and coalgebra with respect to the same grading.

\medbreak
For example, the dual $H^{*}$  of a finite dimensional color Hopf algebra $H$ is again a
color Hopf algebra.

\smallbreak
Let $H$ be a color Hopf algebra. Then $H$ is a Hopf algebra in $^{\F G}_{\F G}\D$ by Remark \ref{rem:color-categories} (iii).
Let $H_G = H\#\F G$ be the  bosonization of $H$, i.~e. $H\otimes \F G$ with multiplication and
coproduct given by \eqref{1.1} and \eqref{1.3}.
Observe that
\begin{align*}
\mathcal{S}(x\#e_g) &= \beta(g^{-1}|x|^{-1},|x|)\mathcal{S}(x)\#e_{g^{-1}|x|^{-1}}, &
x&\in H, \, g\in G.
\end{align*}
Let $_H\mathcal{CM}$ be the category of graded modules with the tensor structure defined as follows.
If $V,W\in\, _H\mathcal{CM}$, then $V\otimes W$ has the underlying tensor product of $G$-graded vector spaces
and the action of $H$ given by
\begin{align}\label{1.6}
    h\cdot(v\otimes w)&= \beta(|h^{(2)}|,|v|)h^{(1)}v\otimes h^{(2)}w,&  v &\in V, \quad w\in W, \quad h\in H.
\end{align}
The functor $_H\mathcal{CM}\to _{H_G}\mathcal{M}$ in Remark \ref{rem:1.1} is monoidal.
The full subcategory of finite dimensional graded $H$-modules is rigid; if $V$ is such
a module, $x\in H_g$, $v^{*}\in V_{h}^{*}=\Hom(V_{h^{-1}},\F)$
and $u\in V_{(gh)^{-1}}$, then  $H$ acts on $V^{*}$ by
$$x\cdot v^{*}:V_{(gh)^{-1}}\to \F,\,\,\,\,\,(x\cdot v^{*})(u)=\beta(g,h)v^{*}(\mathcal{S}(x)u).$$

Similarly, $^H\mathcal{CM}$ denotes the tensor category of graded comodules with the following tensor product.
If $V,W\in\, ^H\mathcal{CM}$, then
\begin{align*}
    \lambda(v\otimes w)&=\beta(|v^{(0)}|,|w^{(-1)}|)v^{(-1)}w^{(-1)}\otimes v^{(0)}\otimes w^{(0)},&  v &\in V, \quad w\in W.
\end{align*}
The equivalence $^H\mathcal{CM}\simeq\, ^{H_G}\mathcal{M}$ in Remark \ref{rem:1.2} is monoidal.
The subcategory of finite dimensional objects in $^H\mathcal{CM}$ is rigid. Indeed, given such a comodule $V$
and $f\in V^{*}$, the coaction $\lambda(f)=f^{(-1)}\otimes f^{(0)}$ is  given by
$f^{(-1)}f^{(0)}(v)=\beta(|v^{(-1)}|,|v^{(0)}|)\,\mathcal{S}^{-1}(v^{(-1)})f(v^{(0)})$, for $v\in V$.

\subsection{Triangular semisimple Hopf algebras and Scheunert's trick}

The problem of classifying  color Hopf algebras with commutation factor
is a  particular case of the problem of classifying Hopf algebras in symmetric categories.
To start with this general problem, we may consider the categories of finite-dimensional representations
of triangular semisimple Hopf algebras. The classification of semisimple triangular Hopf algebras over
$\F$ was established in \cite{EtGe,EtGe1}, whose notations and
definitions we use freely.  We next summarize \cite[Theorem 2.1]{EtGe} and \cite[Theorem 3.1]{EtGe1}.

\begin{theorem}\label{th:eg-classif}
Let $(B,R)$ be a triangular semisimple Hopf algebra over $\F$, with Drinfeld element $u$. Set $R_u:=\frac 12(1\otimes1+1\otimes u+u\otimes 1-u\otimes u)$ and $\widetilde{R}:=RR_u$. Then there exist a finite group $L$, a subgroup $K<L$ and a minimal twist $J\in \F K\otimes \F K$ such that $(B,\widetilde{R})$ and $(\F L^J,J_{21}^{-1}J)$ are isomorphic as triangular Hopf algebras. Moreover, the data $(L,K,J)$ is unique up to isomorphisms of groups and gauge equivalence of twist. \qed
\end{theorem}

There is a well-known one-to-one correspondence between bicharacters of $G$ and quasitriangular structures on $\F A$. For, consider the isomorphism of Hopf algebras $\F^{G}\to \F A$ given by
$\delta_g \mapsto v_{g} :=\frac{1}{|G|}\sum\limits_{a\in A}a(g^{-1}) e_a$; then the bicharacter $\beta$ corresponds to the R-matrix
\begin{align*}
R_{\beta}&=\sum\limits_{g,h\in G}\beta(g, h)\,v_g\otimes v_h.
\end{align*}
Under this identification, the Drinfeld element of $(\F A,R_{\beta})$ is  $u \in \F^{G}$ given by $u(g)=\beta(g,g^{-1})$, $g\in G$.
Also commutation factors of $G$ correspond to triangular structures on $\F A$.
\emph{Assume until the end of this Subsection  that $\beta$ is a commutation factor of $G$}.
Then the Drinfeld element $u \in A \simeq G(\F^{G})$ is an involution, hence $u(g)=\beta(g,g)$, for  $g\in G$, and
\begin{align*} && (R_{\beta})_u &=\sum\limits_{g,h\in G}\kappa(g,h)\,v_g\otimes v_h,
\\ &\text{where} &
\kappa(g,h) &=
\begin{cases}
    1, & \text{if}\,\,\,u(g)=1\,\,\text{or}\,\,u(h)=1,\\
    -1,& \text{if}\,\,\,u(g)=u(h)=-1.
\end{cases}
\end{align*}
Consequently the $R$-matrix $\widetilde{R}_{\beta}$ as in Theorem \ref{th:eg-classif} has the explicit form
$$\widetilde{R}_{\beta}=\frac{1}{|G|^2}\sum\limits_{g,h\in G}\beta\kappa(g,h)\,v_g\otimes v_h.$$
By definition of $\kappa$, $\beta\kappa:G\times G\to \F$ is a commutation factor on $G$. Set
$$ K:= (\ker{\overline{\chi}})^{\perp}=\{ a\in A: a|_{\ker{\overline{\chi}}}=\epsilon_A\},$$
where $\overline{\chi}:G\to A$ is defined by $\overline{\chi}_g(h)=\beta\kappa(g,h)$, $g,h\in G$.
$K$ is identified with the group of characters of $G'=G/\ker{\overline{\chi}}$, and $\beta\kappa$
induces a non-degenerate bicharacter $\beta'$ on $G'$ with $\beta'(\overline{g},\overline{g})=1$ for all $\overline{g}\in G'$. Moreover
\begin{align*}
\widetilde{R}_{\beta}&=\sum\limits_{\overline{g},\overline{h}\in G'}\beta'(\overline{g},\overline{h})\,v_{\overline{g}}
\otimes v_{\overline{h}}, & \text{where }\  v_{\overline{g}} & =\frac{1}{|G'|}\sum\limits_{a\in K} a(\overline{g}^{-1})e_a.
\end{align*}
By Scheunert's trick \cite{Sc0},
there is a normalized 2-cocycle $\gamma'$ on $G'$ such that
\begin{align}\label{eq:trick}
\beta' (\overline{g},\overline{h})&= \frac{\gamma'(\overline{g},\overline{h})}{\gamma'(\overline{h},\overline{g})}, & \text{for all } \overline{g},\overline{h}\in G'.
\end{align}
A well-known and straightforward calculation shows that
$$J_{\gamma'}=\sum\limits_{\overline{g},\overline{h}\in G'}\gamma'(\overline{g},\overline{h})\,v_{\overline{g}}
\otimes v_{\overline{h}} $$
is a minimal twist for $\F K$ such that $\widetilde{R}_\beta=(J_{\gamma'})^{-1}_{21}J_{\gamma'}$;
but $A$ being abelian, $\F A=\F A^{J_{\gamma'}}$ as Hopf algebras. Thus we have:

\begin{theorem} The data associated to the triangular semisimple Hopf algebra $(\F A,R_{\beta})$ in Theorem
\ref{th:eg-classif}  is $(A,K,J_{\gamma'})$.
\qed
\end{theorem}

Since $u^2=1$, it can be thought of as group morphism $u:G\to \mathbb{Z}_2$.
Let $\varsigma$ be the commutation factor on $\Z/2$ such that
$\Vect^{\Z/2}_{\varsigma}$ is the category of supervector spaces; note $\kappa = \varsigma(u\times u)$. Consider the normalized 2-cocycle on $G$ given by $\gamma:=\gamma'(p\times p)$, where $p:G\to G'$ is the canonical projection.
Then $u$ and $\gamma$ give rise to a braided tensor functor
$(F, \phi): \Vect^{G}_\beta \to \Vect^{\Z/2}_{\varsigma}$ by \eqref{eq:trick},
see Remark \ref{rem:color-categories} (i).
Thus, every color Hopf algebra in $\Vect^{G}_{\beta}$ becomes a Hopf superalgebra via $(F, \phi)$.

\smallbreak
In fact, the cocycle $\gamma$ induces  an equivalence of  braided tensor categories
$(F, \phi): \Vect^{G}_{\beta\kappa} \to \Vect^{G}_{\varepsilon}$, where $\varepsilon$ is the trivial bicharacter, see
Remark \ref{rem:color-categories} (i). (This reflects that $(\F A^{J_{\gamma}}, (J_\gamma)^{-1}_{21}J_\gamma)$ is the twist of $(\F A, 1 \otimes 1)$). Hence

\begin{itemize}
	\item If $u = 1$, then $(F, \phi): \Vect^{G}_{\beta} \to \Vect^{G}_{\varepsilon}$ is an equivalence of  braided tensor categories, hence
	any $(G, \beta)$-color Hopf algebra gives rise to a $G$-graded Hopf algebra and vice versa.
	\item If $u \neq 1$,
	any $(G, \beta)$-color Hopf algebra gives rise to a $G$-graded Hopf superalgebra and vice versa.
\end{itemize}
Compare e.~g. with \cite[Section 2]{BD}. By these reasons, and taking into account that in Heckenberger's list \cite{He} there are
Nichols algebras that arise neither from usual Lie algebras nor from Lie superalgebras-- see the discussion in \cite{Ang}-- we
are led to consider bicharacters that are not commutation factors.

\section{Categories of modules over color Hopf algebras}\label{section:color categories}

As always, $G$ is a finite abelian group, $A = \Hom(G, \F^{\times})$,
$\beta$ is a bicharacter on $G$ and $\chi: G \to A$ is  given by \eqref{equation0}. We fix a color Hopf algebra $H$.
The  statements below are basically adaptations of  analogous statements for Hopf algebras (and particular cases
of statements for braided Hopf algebras in braided tensor categories) and we omit most of the proofs.

\subsection{Color Hopf bimodules and Yetter-Drinfeld color modules}

A {\it (left) color Hopf module} over $H$ is an object $V\in \Vect^G$ which
is simultaneously a (left) graded module over $H$ and a (left) graded comodule over $H$ such that the
coaction $\lambda:V\to H\otimes V$ is a morphism of graded modules over $H$; here we use the braiding for the
action on $H\otimes V$, cf. \eqref{1.6}. If $U\in \Vect^G$,
then $H\otimes U$ is a  color Hopf module with the regular left action and coaction. If $V$ is a graded
comodule over $H$, then $V^{\co H}:= \{v\in V\,:\,\lambda(v)=1\otimes v\}$.
The Fundamental theorem of color Hopf modules, a particular case of the analogous theorem
in braided tensor categories \cite{Ta2}, says:
\vu

\begin{theorem}\label{th:2.1} If $V$ is a color Hopf module over $H$ then the multiplication  $H\otimes V^{\co H}\to V$ is
an isomorphism of color Hopf modules.\qed
\end{theorem}

A {\it color Hopf bimodule} over $H$ is an object $V\in \Vect^G$ which is simultaneously a graded bimodule
over $H$ and a graded bicomodule over $H$, and such that the coactions $\lambda:V\to H\otimes V$ and
$\rho:V\to V\otimes H$ are morphisms of graded bimodules over $H$ (for what we need again the braiding).
The category $_H^H\mathcal{CM}_H^H$ of
color Hopf bimodules over $H$ is tensor one with the tensor product $\otimes_H$.
A {\it left Yetter-Drinfeld color module} over $H$ is an object $V\in \Vect^G$ which is simultaneously a left
color module over $H$ and a left color comodule over $H$ with compatibility
$$\lambda(h\cdot v)=\beta(|h^{(2)}||h^{(3)}|,|v^{(-1)}|)\beta(|h^{(2)}|,|h^{(3)}|)h^{(1)}v^{(-1)}
\mathcal{S}(h^{(3)})\otimes h^{(2)}\cdot v^{(0)}
$$
for all $v \in V$ and $h\in H$. As usual, these notions are equivalent.

\begin{prop}\label{prop:yetter} The category $_H^H\mathcal{YDC}$ of left Yetter-Drinfeld color modules is tensor equivalent to $_H^H\mathcal{CM}_H^H$.
\end{prop}

\pf We have functors:
$ _H^H\mathcal{CM}_H^H \ni M\rightsquigarrow V=M^{\co H}\in _H^H\mathcal{YDC}$,
with action $h\cdot v=\beta(|h^{(2)}|,v)h^{(1)}v\mathcal{S}(h^{(2)})$ and coaction $\lambda$; and
$ _H^H\mathcal{YDC} \ni V\rightsquigarrow M=V\otimes H$,
with actions and coactions
\begin{align*}
x\cdot (v\otimes h)&=\beta(|x^{(2)}|,|v|)x^{(1)}v\otimes x^{(2)}h,& (v\otimes h)\cdot y&=v\otimes hy, \\
\lambda(v\otimes h)=\beta(|v^{(0)}|,&|h^{(1)}|)v^{(-1)}h^{(1)}\otimes v^{(0)}\otimes h^{(2)}, &
\rho &= \text{id}\otimes \Delta.
\end{align*}
These functors are tensor equivalent, inverse of each other. \epf

The tensor category $^H_H\mathcal{YDC}$ is braided; explicitly, if $V,W\in ^H_H\mathcal{YDC}$, then the braiding $c_{V,W}:V\otimes W\to W\otimes V$ is
\begin{align}\label{1.7}
&c_{V,W}(v\otimes w)=\beta(|v^{(0)}|,|w|)\,v^{(-1)}w\otimes v^{(0)},& &v\in V,& &w\in W.&
\end{align}

Indeed, $c_{V,W}$ is a natural isomorphism in $^H_H\mathcal{YDC}$ with inverse $$c^{-1}_{V,W}(w\otimes v)=\beta(|v^{(0)}|,|v^{(-1)}w|)^{-1}\,v^{(0)}\otimes\mathcal{S}^{-1}(v^{(-1)})w,\,\,\,v\in V,\,\,w\in W.$$
Thus  $_H^H\mathcal{CM}_H^H$ is braided, with braiding $c_{M,N}:M\otimes_H N\to N\otimes_H M$ given by
$$m\otimes n\mapsto\beta(|m^{(0)}||m^{(-1)}|,|n^{(0)}||n^{(1)}|)\,m^{(-2)}n^{(0)}\mathcal{S}(n^{(1)})\mathcal{S}(m^{(-1)})\otimes m^{(0)}n^{(2)},$$
for all $m\in M$, $n\in N$.

\begin{remark}\label{rem:2.3}
 Let $H$ be a color Hopf algebra and $H_G=H\#\F G$. There exists a full embedding of braided tensor categories $\mathfrak{i}:\,^H_H\mathcal{YDC}\hookrightarrow\, _{H_G}^{H_G}\mathcal{YD}$; indeed, any
 $V\in \,^H_H\mathcal{YDC}$ becomes a Yetter-Drinfeld module over $H_G$ with the structures \eqref{1.2} and \eqref{1.4}.
\end{remark}

\begin{exa}\label{trivial:grading} Let $H$ be a color Hopf algebra with trivial grading, that is, $H$ is a usual Hopf algebra. If $V\in ^H_H\mathcal{YD}$ and $g\in G$, then $V[g]:=V\in ^H_H\mathcal{YDC}$, where all elements of $V[g]$ have degree $g$. Clearly if $V$ is irreducible in $^H_H\mathcal{YD}$, then $V[g]$ is irreducible in $^H_H\mathcal{YDC}$. Reciprocally, if $W=\oplus_{g\in G}W_g$ is an irreducible object in $^H_H\mathcal{YDC}$, then there exists $g\in G$ such that $W=W_g$. Hence, $W=U[g]$ with $U:=W$ an irreducible object in $^H_H\mathcal{YD}$.
\end{exa}

\subsection{Bosonization of color Hopf algebras}

Since $^H_H\mathcal{YDC}$ is a braided category, we may consider Hopf algebras in this category.
The equivalence between Hopf algebras in $^H_H\mathcal{YDC}$ and color Hopf algebras with split projection to $H$
goes in the same way as in \cite{Radford, Majid}.
Namely, given  a Hopf algebra $R$ in $^H_H\mathcal{YDC}$,
the color Hopf algebra $R\#H$  is $R\otimes H\in \Vect^{G}$ with
\begin{align*}
(a\#h)(b\#f)&= \beta(|h^{(2)}|,|b|)a(h^{(1)}\cdot b)\#h^{(2)}f, \\
\Delta(a\#h) &=\beta(|(a^{(2)})^{(0)}|,|h^{(1)}|)a^{(1)}\#(a^{(2)})^{(-1)}h^{(1)}\otimes (a^{(2)})^0\#h^{(2)},  \\
1_{R\#H} &= 1_R\#1_H, \qquad \varepsilon(a\#h)=\varepsilon_R(a)\varepsilon_H(h),\\
\mathcal{S}(a\#h) &=\beta(|a^{(0)}|,|h|)(1_R\#\mathcal{S}_H(a^{(-1)}h))(\mathcal{S}_R(a^{(0)})\#1_H),
\end{align*}
for all $a,b\in R$ and $h,f\in H$.
By Remark \ref{rem:2.3}, $\mathfrak{i}(R)$ is a Hopf algebra in $_{H_G}^{H_G}\mathcal{YD}$.
It is straightforward to verify that
\begin{equation}\label{1.10}
(R\#H)_G\simeq \mathfrak{i}(R)\#H_G.
\end{equation}

Conversely, let $\iota:H\hookrightarrow L$ and $\pi:L\twoheadrightarrow H$ be morphisms of color Hopf algebras such
that $\pi\circ\iota=\id_H$. The subalgebra of coinvariants $$R:=L^{\co H}=\{x\in L\,:\,(\id\otimes\pi)\Delta(x)=x\otimes 1\}$$
is a Hopf algebra in the category $^H_H\mathcal{YDC}$ with the following structures
\begin{align*}
h.r &= \beta(|h^{(2)}|,|r|)\iota(h^{(1)})r\iota(\mathcal{S}(h^{(2)})),& \lambda &=(\pi\otimes id)\Delta, \\
\Delta_R(r) &= r^{(1)}\iota\pi\mathcal{S}(r^{(2)})\otimes r^{(3)}, \quad \varepsilon_R(r)=\varepsilon_L(r),
&\mathcal{S}_R(r) &= \iota\pi(r^{(1)})\mathcal{S}_L(r^{(2)}).
\end{align*}
The color Hopf algebras $L$ and $R\#H$ are isomorphic.

Note that $\iota$ and $\pi$ induce morphisms of Hopf algebras $\iota_G:H_G\hookrightarrow L_G$
and $\pi_G:L_G\twoheadrightarrow H_G$ such that $\pi_G\circ\iota_G=\id_{H_G}$. Hence $\mathfrak{i}(R)$
coincides with the subalgebra of coinvariants $(L_G)^{\co (H_G)}$, see {\it e. g.}  \cite[Lemma 3.1]{AHS}.

\subsection{Color Nichols algebras}

The structure of a Nichols algebra appeared first in \cite{Nichols}. See \cite{AS3} for a detailed introduction
to this construction.

\begin{prop}\label{prop:nichols}  Let $V\in {}^H_H\mathcal{YDC}$. Then there is a unique (up to isomorphisms)
 graded Hopf algebra $\mathfrak{B}(V)=\oplus_{n\in \mathbb{N}_0}\mathfrak{B}^n(V)$ in $^H_H\mathcal{YDC}$
such that:

\begin{enumerate}\renewcommand{\theenumi}{\roman{enumi}}   \renewcommand{\labelenumi}{(\theenumi)}
\item $\mathfrak{B}^0(V)\simeq \F$,

\item  $V\simeq \mathfrak{B}^1(V)=\mathcal{P}(\mathfrak{B}(V))$\,\,\,(the space of primitive elements),

\item  $\mathfrak{B}^1(V)$ generates the algebra $\mathfrak{B}(V)$. \qed
\end{enumerate}
\end{prop}

We have a similar explicit description of $\mathfrak{B}(V)$  as in \cite[1.7]{AAY}.
Then, by Remark \ref{rem:2.3}, it follows that
\begin{equation}\label{1.11}
\mathfrak{B}(\mathfrak{i}(V))=\mathfrak{i}(\mathfrak{B}(V)).
\end{equation}

\subsection{The lifting method for color Hopf algebras}

Let $H$ be a color Hopf algebra. The coradical $H_0$ of $H$ is a graded subcoalgebra of $H$.
Assume that $H_0$ is also a subalgebra of $H$.
Because of this assumption on $H$, the coradical filtration $\{H_n\}$ is also a filtration of algebras and the associated
graded algebra $\text{gr}\,H$ is a color Hopf algebra. Since the homogeneous projection
$\pi:\text{gr}\,H\longrightarrow H_0$ is a morphism of color Hopf algebras that splits the inclusion map,
the algebra of left coinvariants $R=(\text{gr}\,H)^{\co \pi}$ is a
braided Hopf algebra in $^{H_0}_{H_0}\mathcal{YD}$
and $\text{gr}\,H\simeq R\#H_0$; $R=\oplus_{n\geq 0}R(n)$ and
$R(0)=\F1$, $R(1)=\mathcal{P}(R)$.
Consequently the subalgebra generated by $R(1)$ is isomorphic to the Nichols algebra $\mathcal{B}(R(1))$.
That is, we can adapt the lifting method \cite{AS3} to the setting of color Hopf algebras.
However, the classification problems of color Hopf algebras whose coradical is a subalgebra reduce
to classification problems of Hopf algebras with the same property.

\begin{prop}\label{prop:nicholsdim0} Let $H$ be a color Hopf algebra.

(i)  The coradical filtration of $H_G$ is $\{H_n\#\F G\}$.

(ii) $H_0$ is a color Hopf subalgebra of $H$ if and only if $(H_G)_0$ is a Hopf subalgebra of $H_G$. In this case,
\begin{align*}
\text{gr}\,(H_G)\simeq (\text{gr}\, H)_G\simeq (R\#H_0)_G\simeq \mathfrak{i}(R)\#(H_0)_G.\
\end{align*}

(iii) The assignement $H\rightsquigarrow H_G$ provides a one-to-one correspondence between
finite-dimensional strictly pointed color Hopf algebras $H$ and finite-dimensional
pointed Hopf algebras $\Hc$ such that:
	
\begin{itemize}
	\item $G(\Hc)=G(H)\times G$.
	\item  The inclusion $\F G\hookrightarrow \Hc$ has a Hopf algebra retraction
 $\pi:\Hc\longrightarrow \F G$.
	\item As Yetter-Drinfeld module, $\Hc^{\co \pi}$ comes from the grading as in \eqref{eq:explicit-tensorfunctor}.
\end{itemize}
\qed
\end{prop}

\begin{exa}\label{exa:decomp} Let $\Gamma$ be a finite abelian group. Then any finite-dimensional
$V$ in $^{\F\Gamma}_{\F\Gamma}\mathcal{YD}$ has a decomposition in direct sum of sub-objects
of the form $V_{u}^{\zeta}$, $u\in \Gamma$, $\zeta\in \widehat{\Gamma}$, where $u$ determines the
coaction and $\zeta$ the action. Consider $\F\Gamma$ as a color Hopf algebra with trivial grading so
that $^{\F\Gamma}_{\F\Gamma}\mathcal{YDC}$ is a semisimple category. By
Example \ref{trivial:grading}, any finite-dimensional $V\in\, ^{\F\Gamma}_{\F\Gamma}\mathcal{YDC}$
has a basis $x_1,\ldots,x_{\theta}$ such that $x_j\in V_{u_j}^{\zeta_j}[t_j]$, $u_j\in \Gamma$,
$\zeta_j\in\widehat{\Gamma}$ and $t_j\in G$.
\end{exa}

\begin{prop}\label{prop:nicholsdim} Let $u_j\in \Gamma$, $\zeta_j\in\widehat{\Gamma}$,
$t_j\in G$, $1\le j\le \theta$. Set
\begin{align*}
\widetilde{q}_{ij} &= \zeta_j(u_i), & q_{ij} &= \beta(t_i,t_j)\widetilde{q}_{ij}, & 1&\leq i,j\leq \theta.
\end{align*}
Let $V\in {}^{\F\Gamma}_{\F\Gamma}\mathcal{YDC}$ with basis $x_1,\ldots,x_{\theta}$ such that
$x_j\in V_{u_j}^{\zeta_j}[t_j]$. Then the color Nichols algebra $\mathfrak{B}(V)$ has finite dimension
if and only if the connected components of the generalized Dynkin diagram corresponding to the
matrix $(q_{ij})_{1\leq i,j\leq r}$ belong to the list in \cite{He}.
\end{prop}

\pf It follows directly from (\ref{1.7}) and (\ref{1.11}).
\epf

\section{Examples related with quantum doubles}\label{section:drinfeld doubles}

In this section we construct examples of color Hopf algebras from Nichols algebras of diagonal type.
We start by a general construction of quantum doubles extending that of \cite{H-isom}, where
quantum doubles of bosonizations of Nichols algebras by $\F\zt$ are considered.
We describe in Subsection \ref{GC} quantum doubles of realizations over other abelian groups.
A similar construction was studied in \cite{Be} for finite abelian groups; see also \cite{ARS, RaS}.
We show how to give colors to these
quantum doubles, that can be thought of as generalizations of quantized enveloping (super)algebras. Then we show
examples where the color is \emph{consistent} with the associated Weyl groupoid.

\subsection{General construction}\label{GC}

We refer to \cite{KlS, DT} for the definition of skew-Hopf pairing between Hopf algebras $B$ and $B'$ and the associated Hopf algebra
$B \bowtie B'$. Let $\I = \{1, \dots, \theta\}$. We consider a datum
$$\E = \big(\q, \Gamma,\Lambda, (g_i)_{i \in \I},
(\gamma_i)_{i \in \I},(h_i)_{i \in \I}, (\lambda_i)_{i \in \I}, \overline{\mu}  \big),$$ where
\begin{itemize}
  \item  $\q = (q_{ij})_{i, j \in \I}$ has entries in $\F^{\times}$,
  \item $\Gamma,\Lambda$ are two abelian groups,
  \item  $g_i\in\Gamma$, $\gamma_j\in\widehat{\Gamma}$ are such that $\gamma_j(g_i)= q_{ij}$,
  \item  $h_i\in\Lambda$, $\lambda_j\in\widehat{\Lambda}$ are such that $\lambda_j(h_i)=q_{ji}$,
  \item $\overline{\mu}:\Gamma\times\Lambda\to\F^\times$ is a $\Z$-bilinear form such that $\overline{\mu}(g_i,h_j)=q_{ij}$.
\end{itemize}
Let $\overline{\mu}:\F\Gamma\otimes\F\Lambda\to\F$ be the skew-Hopf pairing associated to $\overline{\mu}$.

\smallbreak
For example, \cite[Section 4]{H-isom} deals with the case $\Gamma=\Lambda=\zt$, while we are interested in the
case $\Gamma=\zt\times G$, $\Lambda=\zt\times A$, $G$ our finite abelian group. Another interesting
example is:  $\Gamma$ a finite abelian group, $\Lambda=\widehat{\Gamma}$, where $h_i=\gamma_i$, $\lambda_i=g_i$ under the
canonical identification of $\Gamma$ with the group of characters of $\widehat{\Gamma}$, and $\overline{\mu}$
is the evaluation.

\smallbreak
First we attach to the datum $\E$ two Yetter-Drinfeld modules
\begin{align}\label{eq:def-V}
V&\in\ydg & & \text{with a fixed basis} & &E_1,\ldots,E_\theta, & E_i&\in V^{\gamma_i}_{g_i},
\\ \label{eq:def-W}
W &\in\ydl & &\text{with a fixed basis} & &F_1,\ldots,F_\theta, & F_i &\in W^{\lambda_i}_{h_i}.
\end{align}

Following \cite{H-isom}, we consider the Hopf algebras
\begin{itemize}
  \item $\A=T(V)\#\F\Gamma$, i.~e. the algebra over $\F \Gamma$ generated by elements  $E_i$, $i\in \I$,
with relations $gE_i=\gamma_i(g) E_ig$, and coproduct  determined by:
\begin{align*}
\Delta(g) &= g\otimes g, & g&\in \Gamma, & \Delta(E_i) &= E_i\otimes 1+g_i\otimes E_i, & i&\in \I.
\end{align*}

  \item $\A'=(T(W)\#\F\Lambda)^{\cop}$, the algebra over $\F \Lambda$ generated by elements
$F_i$, $i\in \I$, and relations $hF_i=\lambda_i(h) F_ih$. The coproduct is given by:
\begin{align*}
\Delta(h) &= h\otimes h, & h&\in \Lambda, & \Delta(F_i) &= F_i\otimes h_i + 1\otimes F_i, & i&\in \I.
\end{align*}
\end{itemize}

\begin{prop}\label{prop:pairing TV}
There exists a unique skew-Hopf pairing $\mu:\A\otimes\A'\to\F$ such that, for all $i, j \in \I$, $g\in\Gamma$,
$h\in\Lambda$, $E\in T(V)$ and $F\in T(W)$,
\begin{align}\label{for:cg}
\mu|_{\F\Gamma\otimes\F\Lambda} &= \overline{\mu},& \mu(E_i,F_j) &= -\delta_{i,j},& \mu(E_i,h) &= \mu(g,F_i) = 0,
\\ \label{for:cg2}
&& \mu(Eg,Fh) &= \mu(E,F)\mu(g,h). &&
\end{align}
\end{prop}
\pf
Analogous to that of \cite[Proposition 4.3]{H-isom}. In fact we can define an algebra homomorphism $\Upsilon:\A\to(\A')^\ast$,
by sending $g\mapsto\overline{g}$, $E_i\mapsto\overline{E_i}$, the linear maps determined by the following equations
$$ \overline{g}(Fh)=\epsilon(F)\overline{\mu}(g,h), \qquad \overline{E_i}(F_{j_1}\cdots F_{j_n}h)=-\delta_{n,1}\delta_{j_1,i}.$$
The map $\Upsilon$ is well-defined; indeed, $ \overline{g}\overline{g'}= \overline{gg'}$, $\overline{g}\overline{E_i}
= \gamma_i(g) \overline{E_i}\overline{g}$, for all $g,g'\in\Gamma$, $1\leq i\leq\theta$.
Therefore the bilinear form satisfying \eqref{for:cg} and \eqref{for:cg2} is defined by $\mu(a,b):=\Upsilon(a)(b)$.
The rest of the proof follows as in loc. cit.
\epf

Let $\B=\B(V)\#\F\Gamma$ and $\B'=(\B(W)\#\F\Lambda)^{\cop}$. There exist canonical surjective Hopf algebra
maps $\pi_V:\A\to\B$, $\pi_W:\A'\to\B'$. If $\mathcal{I}(V)$, $\mathcal{I}(W)$ denote the ideals of $T(V)$, $T(W)$
defining the corresponding Nichols algebras, then
$$\ker \pi_V=\mathcal{I}(V)\Gamma, \qquad \ker \pi_W=\mathcal{I}(W)\Lambda.$$
We prove now \cite[Theorem 5.8]{H-isom} in our general setting; cf. also \cite{ARS, Be}.

\begin{theorem}\label{thm:pairing nichols}
The skew-Hopf pairing $\mu$ of Proposition \ref{prop:pairing TV} induces a skew-Hopf pairing
$\mu:\B\otimes\B'\to\F$, whose restriction to $\B(V)\otimes\B(W)$ is non-degenerate; $\mu$ is non-degenerate if and only if $\overline{\mu}$ is non-degenerate.
\end{theorem}

\pf
Set $\Vc^+=T(V)\#\zt$ (and call $(K_i)_{i\in \I}$ the canonical basis of this $\zt$);
$\Vc^-=T(W)\#\zt$ (and call $(L_i)_{i\in \I}$ the canonical basis of this second copy of $\zt$), $\nu:\Vc^+\otimes\Vc^-\to \F$ the bilinear form as in Proposition \ref{prop:pairing TV} but with respect with these copies of $\zt$, cf. \cite[Proposition 4.3]{H-isom}.
There exist Hopf algebra maps $\alpha_V:\Vc^+\to\A$, $\alpha_W:\Vc^-\to\A'$, determined by the conditions
$\alpha_V(K_i) = g_i$, $\alpha_V|_{T(V)} = \id_{T(V)}$,
$\alpha_W(L_i)=h_i$, $\alpha_W|_{T(W)}=\id_{T(W)}$,
$i\in \I$. Then
\begin{equation}\label{eq:pairings heck y generalizado}
\nu(x,y)=\mu(\alpha_V(x),\alpha_W(y)), \qquad \mbox{for all }x\in\Vc^+,\, y\in\Vc_-.
\end{equation}
We claim that $\mu(\mathcal{I}(V)\Gamma,\A')=0$, $\mu(\A,\mathcal{I}(W)\Lambda)=0$. To prove this, we use \eqref{for:cg2}
and argue as in the proof of \cite[Theorem 5.8]{H-isom}.
Therefore we induce a bilinear form $\mu:\B\otimes\B'\to\F$,
which is a skew-Hopf pairing because $\pi_V$, $\pi_W$ are Hopf algebra morphisms.
The restriction of the pairing to $\B(V)\otimes\B(W)$ is non-degenerate because of
\eqref{eq:pairings heck y generalizado}
and again the proof of \cite[Theorem 5.8]{H-isom}.
The last claim follows from the preceding and \eqref{for:cg2}.
\epf

\begin{definition}\label{def:doble drinfeld nichols}
Let $D(\E)$ be the Hopf algebra obtained from  $\B\otimes\B'$ twisting by the cocycle associated to $\mu$ \cite{DT, KlS}.
It is presented by generators $g\in\Gamma$, $h\in\Lambda$, $E_i, F_i$, $i\in \I$,
the relations in $\mathcal{I}(V)$, $\mathcal{I}(W)$,
those of $\Gamma\times\Lambda$,  and
\begin{align*}
gE_i &= \gamma_i(g) E_ig, &   hE_i &= \mu(g_i,h)^{-1} E_ih,
\\ hF_i &= \lambda_i(h) F_ih,  &  gF_i &= \mu(g,h_i)^{-1} F_ig,
\\ E_iF_j-F_jE_i &= \delta_{ij}(g_i-h_i), & i,j&\in \I, g\in\Gamma, h\in\Lambda.
\end{align*}
\end{definition}

\subsection{Coloring quantum doubles}
We now define a color version of the Hopf algebra $D(\E)$ above, for the
pair $(G, \beta)$; we assume that $\beta$ is non-degenerate in this Subsection.
Hence the maps $\chi, \chi^o: G \to A$ as in  \eqref{equation0}
are isomorphisms of groups. To start with, we consider the following datum:
\begin{itemize}
  \item A matrix  $\q = (q_{ij})_{i, j \in \I}\in (\F^{\times})^{\I\times \I}$, such that $q_{ii}\neq 1$, for all $i \in \I$.
  \item $t_i\in G$ for all $i \in \I$.
\item $\Gamma=\Go\times G$, $\Go$ free abelian with basis $K_1,\ldots,K_{\theta}$.
  \item $\Lambda=\Lo\times A$, $\Lo$ free abelian with basis $L_1,\ldots,L_{\theta}$.
\end{itemize}
Set $\xi_i=\chi_{t_i}\in A$, $\widetilde{q}_{ij} = \beta(t_i,t_j)^{-1}q_{ij}
= \xi_j(t_i)^{-1}q_{ij}$, $i,j\in \I$. 	
Let $\E = \big(\q, \Gamma,\Lambda, (g_i)_{i \in \I},
(\gamma_i)_{i \in \I},(h_i)_{i \in \I}, (\lambda_i)_{i \in \I}, \overline{\mu}  \big)$, where

\begin{itemize}	\renewcommand{\labelitemi}{$\circ$}
 \item  $g_i = (K_i, t_i)\in\Gamma$.
\item $\gamma_j = (\zeta_j, \xi_j)\in\widehat{\Gamma}$, where
$\zeta_j\in \widehat{\Go}$ is given by $\zeta_j(K_i)= \widetilde{q}_{ij}$, $j\in \I$.
 \item  $h_i = (L_i, \xi_i)\in\Lambda$.
\item  $\lambda_j = (\eta_j, t_j)\in\widehat{\Lambda}$
where $\eta_j\in \widehat{\Lo}$ is given by $\eta_j(L_i)= \widetilde{q}_{ji}$,  $i\in \I$.
  \item $\overline{\mu}:\Gamma\times\Lambda\to\F^\times$ is given by
$\overline{\mu}\big((K_i, t), (L_j, \alpha) \big) = \widetilde{q}_{ij} \alpha(t)$.
\end{itemize}

We shall say that $\E$ is a datum, or a datum for $\q$ if emphasis is needed.

\medbreak
Let $D(\E)$ be the Hopf algebra as in Definition \ref{def:doble drinfeld nichols},
presented by generators $E_i$, $F_i$, $K_i^{\pm 1}$, $L_i^{\pm 1}$, $i\in\I$, $t\in G$,
$\xi\in A$ with the relations defining the Nichols algebras generated by $E_i$  and $F_i$,
$i\in\I$; those of $\Gamma\times\Lambda$, and
\begin{align}
\label{rel1}
&K_iE_jK_i^{-1}=\widetilde{q}_{ij}E_j,&   &K_iF_jK_i^{-1}=\widetilde{q}_{ij}^{-1}F_j,& \\ \label{rel2}
&L_iE_jL_i^{-1}=\widetilde{q}_{ji}^{-1}E_j,& &L_iF_jL_i^{-1}=\widetilde{q}_{ji}F_j,&  \\ \label{rel3}
&tE_jt^{-1}=\xi_j(t)E_j = \beta(t,t_j)E_j,&  &tF_jt^{-1}=\xi_j(t)^{-1}F_j,&\\ \label{rel4}
&\xi E_j\xi^{-1}=\xi(t_j)^{-1}E_j,& &\xi F_j\xi^{-1}=\xi(t_j)F_j,& \\ \label{rel5}
&E_iF_j-F_jE_i=\delta_{ij}(t_iK_i-\xi_i L_i),& &i, j\in \I, t\in G, \xi\in A.&
\end{align}

Here the elements of
$\Gamma\times\Lambda$ are group-like and
\begin{align}\label{comul}
&\Delta(E_i)=E_i\otimes 1+K_it_i\otimes E_i,& & \Delta(F_i)=F_i\otimes L_i\xi_i+1\otimes F_i.&
\end{align}

\medskip

Given $t\in G$, $t\chi_t^o$ is a central element of $D(\E)$. Let $I(\E)$ be the Hopf ideal of $U(\E)$
generated by $\{t\chi_t^o-1,\,\, t\in G\}$ and $U(\E):= D(\E)/I(\E)$ the
corresponding quotient Hopf algebra; it has a presentation by generators $E_i$, $F_i$, $K_i^{\pm 1}$, $L_i^{\pm 1}$, $i\in\I$, $t\in G$ and
analogous relations.
Our goal is to obtain color Hopf algebras from a factorization $U(\E)\simeq R(\E)\# \F G$. Therefore we analyze first such possible factorizations.

\begin{prop}\label{prop:quot} There exists a one-to-one correspondence between
\begin{enumerate}\renewcommand{\theenumi}{\roman{enumi}}   \renewcommand{\labelenumi}{(\theenumi)}
\item group homomorphisms $\overline{\pi}:\Go\to G$ and
\item Hopf algebra retractions  $\pi:U(\E)\longrightarrow \F G$
of  $\F G\hookrightarrow U(\E)$.
\end{enumerate}
Indeed, given $\overline{\pi}:\Go\to G$, we define $\pi:U(\E)\to \F G$ by $\pi(t)=t$, $t\in G$, and $\pi(K_i)=\overline{\pi}(K_i)$,
$\pi(L_i)=t_i\overline{\pi}(K_i)\chi^{-1}(\xi_i)$, $\pi(E_i)=\pi(F_i)=0$, $i\in\I$.
\end{prop}

\pf It is easy to see that the extension of $\overline{\pi}:\Go\to \F G$
to  $\pi:U(\E)\to \F G$ as above is a Hopf algebra map.
Conversely, let $\pi:U(\E)\to \F G$ be a Hopf algebra retraction of $\F G \hookrightarrow U(\E)$. Since $\pi(K_i),\pi(L_i)$
are group-likes, $\pi|_{\Go\times\Lo}$ is a group homomorphism.
By \eqref{rel1}, \eqref{rel3} and the hypothesis
$q_{ii}\neq 1$,  $\pi(E_i)=0 =\pi(F_i)$.
Hence $\pi(L_i)=t_i\overline{\pi}(K_i)\chi^{-1}(\xi_i)$ by \eqref{rel5}.
\epf

Let $\pi:U(\E)\to \F G$ be a Hopf algebra map as in  Proposition \ref{prop:quot}
and $R(\E)=U(\E)^{\co \pi}$, a braided Hopf algebra in $\ydG$. Notice that
$R(\E)$ is the subalgebra of
$U(\E)$ generated by $E_i$, $\widetilde{K}_i=K_i\pi(K_i)^{-1}$,
$\widetilde{L}_i=L_i\pi(L_i)^{-1}$, $\widetilde{F}_i=F_it_i^{-1}\pi(K_i)^{-1}$, $i\in\I$.
Indeed the later is contained in the former; but the $\zt$-homogeneous components coincide,
from $U(\E)=R(\E)\#\F G$ and
$ U(\E) \simeq U^+(\chi,\beta)\otimes (\F\mathbb{Z}^{2\theta}\otimes \F G)\otimes U^{-}(\chi,\beta)$.

Next we decide when the braided Hopf algebra $R(\E)$ is color.

\begin{prop}\label{prop:Drincolor} $R(\E)$ is a color Hopf algebra if and only
if $\overline{\pi}$ is trivial.
\end{prop}
\pf
Since $\lambda=(\pi\otimes \text{id})\Delta$ is the coproduct in $R(\E)\# \F G$,
$\lambda(\widetilde{K}_i)=1\otimes \widetilde{K}_i$, $\lambda(\widetilde{L}_i)=1\otimes \widetilde{L}_i$,
$\lambda(E_i)=\overline{\pi}(K_i)t_i\otimes E_i$ and $\lambda(\widetilde{F}_i)=t_i^{-1}\overline{\pi}(K_i)^{-1}
\otimes \widetilde{F}_i$. We have to check that the action is given by \eqref{eq:explicit-tensorfunctor}.
But $t\cdot E_i=\beta(t,t_i)E_i$ and $t\cdot \widetilde{F}_i=\beta(t,t_i^{-1})\widetilde{F}_i$, $t\in \Gamma$.
Since $\beta$ is non-degenerate, $R(\E)$ is a color Hopf algebra if and only if $\overline{\pi}(K_i)=1$.
\epf

Proposition \ref{prop:Drincolor} gives many examples of color Hopf algebras, namely $R(\E)$, by
taking $\overline{\pi}$ trivial.

\smallbreak
An important invariant of Nichols algebras of diagonal type is the Weyl groupoid, see
\cite{He} for more details and \cite{AHS} for a generalization.
Let $\q = (q_{ij})_{i, j \in \I}\in (\F^{\times})^{\I\times \I}$, such that $q_{ii}\neq 1$, for all $i \in \I$, as above. Define $a_{pp}=2$ and for $j\neq p\in\I$
\begin{align}\label{eq:def aij}
 a_{pj}&:=-\min\{ n: (1+q_{pp}+\cdots+q_{pp}^n)(1-q_{pp}^nq_{pj}q_{jp})=0 \}.
\end{align}
Let $p\in\I$ such that all the $a_{pj}$'s are finite, $j\in\I$.
We define
\begin{itemize}
\item reflections $s_p:\Go\to\Go$ and $s_p:\Lo\to \Lo$ by
\begin{align}\label{eq:reflections}
 s_p(K_j)&=K_jK_p^{-a_{pj}}, & s_p(L_j)&=L_jL_p^{-a_{pj}}, & i\in&\I,
\end{align}
\item the $p$-reflected matrix $s_p^\ast \q= (q_{ij}')$,  where
\begin{align}\label{eq:def reflection}
 q_{ij}'&= q_{ij}q_{pj}^{-a_{pi}}q_{ip}^{-a_{pj}}q_{pp}^{a_{pi}a_{pj}}, &  i,j&\in\I.
\end{align}
\end{itemize}

The Weyl equivalence relation
on the matrices $\q$ as above is induced by $\q\sim s_p^\ast \q$,
for each $p$ such that all the $a_{pj}$'s are finite.
Indeed, $s_p^\ast (s_p^\ast \q)=\q$ \cite{H-isom}. We now extend the Weyl equivalence to data $\E$ as above by the following rule.
Let $p\in\I$ such that all the $a_{pj}$'s are finite, $j\in\I$. Let

$$s_p^{\ast}\E = \E' = \big(\q', \Gamma,\Lambda, (g_i')_{i \in \I},
(\gamma_i')_{i \in \I},(h_i')_{i \in \I}, (\lambda_i')_{i \in \I}, \overline{\mu}  \big),$$
be the datum given by

\bigbreak
\begin{itemize}
\item $\q' = s_p^{\ast}\q$.
\item $g^{\prime}_i=g_ig_p^{-a_{pi}}$; thus the basis of $\Go$ for $\E'$ is
$K_1' = s_p(K_1), \dots, K_\theta' = s_p(K_\theta)$, while $t^{\prime}_i= t_it_p^{-a_{pi}}$.
It follows that $\widetilde{q}_{ij}^{\,\,\prime} = \beta(t'_i,t'_j)^{-1}q'_{ij}$, satisfies
\begin{align*}
\widetilde{q}^{\,\,\prime}_{ij} &= \widetilde{q}_{ij}(\widetilde{q}_{pj})^{-a_{pi}}
(\widetilde{q}_{ip})^{-a_{pj}} (\widetilde{q}_{pp})^{a_{pi}a_{pj}}, & i,j&\in \I.
\end{align*}
\item $\gamma^{\prime}_j = (\zeta'_j, \xi'_j)\in\widehat{\Gamma}$, where $\xi'_j = \chi_{t'_j}$
and $\zeta'_j\in \widehat{\Go}$ is given by
$\zeta'_j(K'_i)= \widetilde{q}^{\,\,\prime}_{ij}$, $j\in \I$.
Hence $\xi^{\prime}_j=\xi_j\xi_p^{-a_{pj}}$, $\zeta_j'=\zeta_j\zeta_p^{-a_{pj}}$
and then $\gamma^{\prime}_j=\gamma_j\gamma_p^{-a_{pj}}$.
\item $h_i'=h_ih_p^{-a_{pi}}$; thus the basis of $\Lo$ for $\E'$ is
$L_1' = s_p(L_1), \dots, L_\theta' = s_p(L_\theta)$.
\item $\lambda'_j = (\eta'_j, t'_j)\in\widehat{\Lambda}$
where $\eta'_j\in \widehat{\Lo}$ is given by $\eta'_j(L'_i)= \widetilde{q}^{\,\,\prime}_{ji}$,
$i\in \I$. Hence $\eta_j' =\eta_j\eta_p^{-a_{pj}}$ and $\lambda_j'=\lambda_j\lambda_p^{-a_{pj}}$.
\end{itemize}
Notice that the bilinear form $\overline{\mu}$ from $\E$ satisfies
$\overline{\mu}\big((K'_i, t), (L'_j, \alpha) \big) = \widetilde{q}_{ij}^{\,\,\prime} \alpha(t)$.

\begin{definition} Let $(\q^{(b)})_{b\in B}$ be the Weyl equivalence class of $\q$.
A family of data $\E^{(b)}$ for $\q^{(b)}$, for $b\in B$, is
\emph{consistently colored} when $\E^{(d)} = s_p^{\ast}\E^{(b)}$ for every two adjacent $b, d\in B$ (that is
$\q^{(d)} = s_p^{\ast}\q^{(b)}$ for some $p$).
\end{definition}

The construction above shows that a datum $\E$ for $\q$ determines a consistently colored family of data $\E^{(b)}$ for the Weyl equivalence class of $\q$.

\bigbreak
 Recall from \cite{He} that the generalized Dynkin diagram of $\q$ has set of vertices $\I$,
each vertex labeled with $q_{ii}$, and arrows between $i,j$ only when $q_{ij}q_{ji}\neq 1$, labeled with this scalar.
\smallbreak

We introduce also \emph{color Dynkin diagrams} for the corresponding matrices $\widetilde{\q}$. It is a diagram with vertices in $\I$,
where each vertex $i$ is labeled with the degree $t_i$ on $G$ and the scalar $\widetilde{q}_{ii}$. Now we put the following set of arrows:
\begin{itemize}
 \item an arrow between $i,j\in\I$ when $\widetilde{q}_{ij}\widetilde{q}_{ji}\neq 1$, labeled with this scalar;
 \item an arrow between $i,j\in\I$ when $\widetilde{q}_{ij}\widetilde{q}_{ji}= 1$ but $\beta(t_i,t_j)\beta(t_j,t_i)\neq 1$, with no labels;
 \item no arrow between $i,j\in\I$ when $\widetilde{q}_{ij}\widetilde{q}_{ji}=\beta(t_i,t_j)\beta(t_j,t_i)= 1$; in particular, when
$\widetilde{q}_{ij}\widetilde{q}_{ji}= 1$, and either $t_i=e$ or $t_j=e$.
\end{itemize}

We show now examples of consistently colored families of data.

\begin{exa}\label{example: C2xC2} Consider the abelian group $G=C_2 \times C_2=\langle\sigma\rangle\times \langle\nu\rangle$,
and the non-degenerate bicharacter $\beta:G\times G\to \F$ given by:
\begin{align*}
 \beta(\sigma,\sigma)&=\beta(\nu,\nu)=\beta(\sigma,\nu)=-1, & \beta(\nu,\sigma)&=1.
\end{align*}
We fix the following notation for the $G$-grading on colored Dynkin diagrams:
\begin{align*}
&\circ\mbox{ for degree }1, &  &\bullet\mbox{ for degree }\sigma, &  & \otimes\mbox{ for degree }\nu, &  & \odot\mbox{ for degree }\sigma\nu.
\end{align*}
Fix $q\in\F^{\times}$, $q\neq\pm 1$, and a matrix $(q_{ij})_{1\leq i,j\leq 4}$ associated to the generalized Dynkin diagram
of \cite[Table 3, Row 14]{He}:
\begin{align}\label{diag:dynkin}
\xymatrix{\circ^{q}\ar@{-}[r]^{q^{-1}}&\circ^{-1} \ar@{-}[r]^{-1}&\circ^{-1} \ar@{-}[r]^{-q}&\circ^{-q^{-1}}}
\end{align}
Set $t_1=e$, $t_2=t_4=\sigma$, $t_3=\nu$.
The matrix $(\widetilde{q}_{ij})_{1\leq i,j\leq 4}$ as in Proposition \ref{prop:nicholsdim} has the colored Dynkin diagram
\begin{align*}
\xymatrix{\circ^{q}\ar@{-}[r]^{q^{-1}}&\bullet^{1} \ar@{-}[r]&\otimes^{1} \ar@{-}[r]^{q}&\bullet^{q^{-1}}}
\end{align*}

In the following table we present the associated generalized and colored Dynkin diagrams of the associated consistently colored family of data.

\vspace{0.4cm}
\noindent
\begin{tabular}{|c|c|}
\hline
 $(q_{ij})=(\beta(t_i,t_j)\widetilde{q}_{ij})$ & $(\widetilde{q}_{ij})$ and $G$-grading \\
\hline
 \xymatrix{ & \circ^{-1} \ar@{-}[rd]^{-1}& &\\\circ^{-1} \ar@{-}[ru]^{q}\ar@{-}[rr]^{-q^{-1}} &  & \circ^{-1}\ar@{-}[r]^{-q} & \circ^{-q^{-1}} } &
 \xymatrix{ & \bullet^{1} \ar@{-}[rd]& &\\\bullet^{1} \ar@{-}[ru]^{q}\ar@{-}[rr]^{q^{-1}} &  & \odot^{1}\ar@{-}[r]^{q} & \bullet^{q^{-1}} }  \\
\hline
 \xymatrix{\circ^{q}\ar@{-}[r]^{q^{-1}} & \circ^{-1} \ar@{-}[r]^{-q} & \circ^{-q^{-1}} \ar@{-}[r]^{-q}& \circ^{-q^{-1}}} &
 \xymatrix{\circ^{q}\ar@{-}[r]^{q^{-1}} & \bullet^{1} \ar@{-}[r]^{q} & \otimes^{q^{-1}} \ar@{-}[r]^{q}& \bullet^{q^{-1}}} \\
\hline
\xymatrix{ &  & \circ^{-1} \ar@{-}[rd]^{q}& \\ \circ^{-q^{-1}}\ar@{-}[r]^{-q} & \circ^{-1} \ar@{-}[ru]^{-1}\ar@{-}[rr]^{-q^{-1}} & & \circ^{-1} } &
\xymatrix{ &  & \otimes^{1} \ar@{-}[rd]^{q}& \\ \otimes^{q^{-1}}\ar@{-}[r]^{q} & \odot^{1} \ar@{-}[ru]^{1}\ar@{-}[rr]^{q^{-1}} & & \otimes^{1} }\\
\hline
\xymatrix{ &  & \circ^{-1} \ar@{-}[rd]^{-q^{-1}}& \\ \circ^{q}\ar@{-}[r]^{q^{-1}} & \circ^{-1} \ar@{-}[ru]^{-1}\ar@{-}[rr]^{q} & & \circ^{-1} } &
\xymatrix{ &  & \otimes^{1} \ar@{-}[rd]^{q^{-1}}& \\ \circ^{q}\ar@{-}[r]^{q^{-1}} & \odot^{1} \ar@{-}[ru]\ar@{-}[rr]^{q} & & \odot^{1} }\\
\hline
\xymatrix{ & \circ^{-1} \ar@{-}[rd]^{-1}& &\\\circ^{-1} \ar@{-}[ru]^{-q^{-1}}\ar@{-}[rr]^{q} &  & \circ^{-1}\ar@{-}[r]^{q^{-1}} & \circ^{q} } &
\xymatrix{ & \bullet^{1} \ar@{-}[rd]& &\\\odot^{1} \ar@{-}[ru]^{q^{-1}}\ar@{-}[rr]^{q} &  & \odot^{1}\ar@{-}[r]^{q^{-1}} & \circ^{q} }\\
\hline
\xymatrix{\circ^{q}\ar@{-}[r]^{q^{-1}} & \circ^{q} \ar@{-}[r]^{q^{-1}} & \circ^{-1} \ar@{-}[r]^{-q}& \circ^{-q^{-1}}} &
\xymatrix{\circ^{q}\ar@{-}[r]^{q^{-1}} & \circ^{q} \ar@{-}[r]^{q^{-1}} & \odot^{1} \ar@{-}[r]^{q}& \bullet^{q^{-1}}}\\
\hline
\end{tabular}
\end{exa}

\begin{exa}\label{nuestro}
Let $\sigma$ be a generator of the cyclic group $G=C_3$ of order $3$, $\omega$ a primitive third root of unity and $q\in\F^\times$, $q\neq 1$,
$\omega$, $\omega^2$.
Let $(\widetilde{q}_{ij})_{1\leq i,j\leq 2}$ be the matrix given by $\widetilde{q}_{ij}=q^{b_{ij}}$, where $b_{11}=0$,
$b_{12}=b_{21}=-1$, $b_{22}=2$.
Let $\beta$ be the bicharacter on $G$ such that $\beta(\sigma,\sigma)=\omega$ and $t_1=\sigma$, $t_2=e$. The matrix $(q_{ij}=\beta(t_i,t_j)\widetilde{q}_{i,j})_{1\leq i,j\leq 2}$ has the generalized Dynkin diagram
\begin{align*}
\xymatrix{\circ^{\omega}\ar@{-}[r]^{q^{-1}}&\circ^{q}}.
\end{align*}
Its Weyl equivalence class includes only one more diagram, namely
\begin{align*}
\xymatrix{\circ^{\omega}\ar@{-}[r]^{q\omega^2}&\circ^{q^{-1}\omega}}.
\end{align*}
Then $t_1'=t_2'=\sigma^2$. If $\circ$, $\bullet$, $\otimes$
denote the vertices of degree $1$, $\sigma$, $\sigma^2$, respectively, the corresponding colored Dynkin diagrams are
\begin{align*}
&\xymatrix{\bullet^{1}\ar@{-}[r]^{q^{-1}}&\circ^{q}}  &
&\xymatrix{\otimes^{1}\ar@{-}[r]^{q}&\otimes^{q^{-1}}}.&
\end{align*}
\end{exa}

\smallbreak

To conclude this subsection, we discuss examples of color Hopf algebras arising from a version of $U(\E)$ but with only one copy of $\zt$.
More precisely, we assume that $(\widetilde{q}_{ij})_{1\leq i,j\leq \theta}$ is symmetric, so
$K_iL_i$ is a central element on $U(\E)$, for all $i \in \I$. Fix $\mathcal{U}(\E):=U(\E)/J$,
where $J$ is the Hopf ideal of $U(\E)$ generated by the set $\{K_iL_i-1\,:\,i \in \I\}$.
The next proposition is an immediate consequence of Propositions \ref{prop:quot} and \ref{prop:Drincolor}.
There is a condition on the elements $t_i$,
not always satisfied as we will see in Example \ref{Yamane}.

\begin{prop}\label{prop:quotcolor} Let $\mathcal{U}(\E)$ be as above.

\noindent \emph{(i)} There exists a one-to-one correspondence between
\begin{enumerate}\renewcommand{\theenumi}{\alph{enumi}}   \renewcommand{\labelenumi}{(\theenumi)}
\item group morphisms $\overline{\pi}:\mathbb{Z}^{\theta}\to G$, such that $\overline{\pi}(K_i^2)=t_i^{-2}$,
\item surjective Hopf algebra morphisms $\pi: \mathcal{U}(\E)\to \F G$ that splits the natural inclusion
$\F G \hookrightarrow \mathcal{U}(\E)$.
\end{enumerate}

\noindent \emph{(ii)} $\mathcal{U}(\E)^{\co \pi}$ is a color Hopf algebra if and only if $t_i^2=1$, for all $1\leq i\leq\theta$.
 \qed
\end{prop}

\begin{exa}\label{Yamane} The next algebra is described in \cite{Y}. It corresponds to our construction in
Proposition \ref{prop:quotcolor} (i) and gives place to a Yetter-Drinfeld Hopf algebra, which is not a color Hopf
algebra.
Let $(\widetilde{q}_{ij})_{1\leq i,j\leq 2}$, $\sigma$ and $\omega$ be as in Example
\ref{nuestro}. The algebra $\mathcal{U}(\E)$ is presented by generators
$E_i$, $F_i$, $K_i^{\pm 1}$, $\sigma$, $i=1,2$,  relations \eqref{rel1}, relations from the group and
\begin{align*}
K_iE_j&=q^{b_{ij}}E_jK_i, & K_iF_j&=q^{-b_{ij}}F_jK_i
\\ \sigma E_i&=\omega^{\delta_{1i}}E_i\sigma,
&\sigma F_i&=\omega^{-\delta_{1i}}F_i\sigma
\\ E_1^3&=(\text{ad}_cE_2)^2E_1=0 & F_1^3&=(\text{ad}_cF_2)^2F_1=0,
\\ E_iF_j-F_jE_i &=\delta_{ij}(-K_i\sigma^{\delta_{1i}}+K_i^{-1}\sigma^{-\delta_{1i}}).& &
\end{align*}
Note that $\pi:\mathcal{U}(\E)\to \F C_3 $, given by $\pi(\sigma)=\sigma$, $\pi(E_i)=\pi(F_i)=0$ and
$\pi(K_i)=\sigma^{-\delta_{1i}}$, is a Hopf algebra morphism that splits the natural inclusion $\F C_3\hookrightarrow \mathcal{U}(\E)$.
In our context $t_1=\sigma$, $t_2=1$, so $t_1^2\neq 1$ and $\mathcal{U}(\E)^{\co \pi}$ is not a color Hopf algebra by Proposition
\ref{prop:quotcolor} (ii).
\end{exa}

\smallbreak

\section{Examples of semisimple color Hopf algebras}\label{section:semisimple color}

Throughout this section, $G$, $A$, $\beta$ and $\chi_g$ are as above, see \eqref{equation0}.
If $L< G$ is a subgroup, then $L^{\bot}:=\{a\in A\,:\,a|_{L}=\varepsilon\}$. Clearly, if $G_1,\ldots,G_r$ are subgroups of $G$,
then $(\cap_{i=1}^r G_i)^{\bot}=\left\langle G_1^{\bot}\cup\ldots\cup G_r^{\bot}\right\rangle$.
We start by a general fact from \cite[Corollary 2]{Schau}; see also \cite[formula (9), page 334]{Du}.

\begin{prop}\label{prop:sch}
Let $H$ and $P$ be Hopf algebras. Assume that $H$ is an object in $_P^P\mathcal{YD}$ for some action and coaction, in such a way that it is an algebra and a coalgebra in $_P^P\mathcal{YD}$. Then the following are equivalent:

\begin{enumerate}\renewcommand{\theenumi}{\roman{enumi}}   \renewcommand{\labelenumi}{(\theenumi)}
\item $H$ is a Hopf algebra in $_P^P\mathcal{YD}$.

\item The braiding $c:H\otimes H\to H\otimes H$, $c(x\otimes y)=x^{(-1)}\cdot y\otimes x^{(0)}$, $x,y\in H$, is the usual flip. \qed
\end{enumerate}
\end{prop}

\begin{remark}\label{rem:gral-hopf-usual-color} Let $H$ be a Hopf algebra.
Then  $H$ is a $G$-graded Hopf algebra if and only if there exists a morphism
$\varphi:A\to \Aut H$ (i.~e. an action of $A$ on $H$ by Hopf algebra automorphisms); clearly
$\Stab_A(H) = \langle \sup H\rangle^{\bot}$.
Hence  $H$ is a $(G,\beta)$-color Hopf algebra if and only if $\beta|_{\,\text{sup\,H}\times \text{sup\,H}}=1$,
if and only if $\chi(\supp H) < \Stab_A(H)$.
\end{remark}

\subsection{Group algebras}

In this subsection, we determine when a group algebra $\F\Gamma$, $\Gamma$ a finite group, is a color Hopf algebra.
By Remark \ref{rem:gral-hopf-usual-color}, $\F\Gamma$ is a $G$-graded Hopf algebra  if and only if there exists
an action by group automorphisms of $A$ on $\Gamma$. Fix such an action; then also $G$ acts on $\Gamma$ via $\chi$.
Clearly $\Stab_A(\F\Gamma) = \Stab_A(\Gamma)$, hence $\beta$ is trivial on $\supp \F\Gamma$ if and only if
$$\supp \F\Gamma < \chi^{-1}(\Stab_A(\Gamma)) = \Stab_G(\Gamma) = \cap_{\eta\in \Gamma} G^{\eta}. $$

\begin{definition}\label{def:3.4} A {\it $(G,\beta)$-color group} is a group $\Gamma$ provided with an action of $A$ by group automorphisms such that $(G^{\eta})^{\bot}< \Stab_A (\Gamma)$, for all $\eta\in \Gamma$.
\end{definition}

\begin{theorem} Let $\Gamma$ be a group provided with an action of $A$ by group automorphisms. Then the following are equivalent:

\begin{enumerate}\renewcommand{\theenumi}{\roman{enumi}}   \renewcommand{\labelenumi}{(\theenumi)}
\item $\F\Gamma$ is a color Hopf algebra,

\item $\beta(g,h)=1$, for all $g,h\in \text{sup}\,\,\F\Gamma$,

\item $\Gamma$ is a $(G,\beta)$-color group,

\item $\F^{\Gamma}$ is a color Hopf algebra.
\end{enumerate}

\end{theorem}

\pf (i) $\Leftrightarrow$ (ii) follows by from Proposition \ref{prop:sch}; (i) $\Leftrightarrow$ (iv) is clear.
(ii) $\Rightarrow$ (iii): Consider $\eta\in \Gamma$ and $a\in (G^{\eta})^{\bot}$. Since  $\text{sup}\,\F\Gamma\subset  G^{\eta}$, $a\in \Stab_A (\Gamma)$.
(iii) $\Rightarrow$ (ii): If $g\in \supp{\F\Gamma}$, then $\chi_g(t)=\chi_t(g)^{-1}=1$, for all $t\in \chi^{-1}(\Stab_A(\Gamma))$. Hence, $\chi_g\in (\cap_{\eta\in\Gamma}G^{\eta})^{\bot}=\left\langle\cup_{\eta\in\Gamma}(G^{\eta})^{\bot}\right\rangle\subset \Stab_A(\Gamma)$, and consequently, $\beta(g,h)=\chi_g(h)=1$, for all $h\in \supp{\F\Gamma}$.\epf

\begin{exa}\label{colorgroup1} Let $G= C_4 = \langle g\rangle$ be the cyclic group of order 4,
${\mathbf i}\in \F$ such that ${\mathbf i}^2=-1$ and $\beta:G\times G \to \F^{\times}$, $\beta(g^i,g^j)=(-1)^{ij}$;
then $A=\langle a\rangle$, where $a(g)={\mathbf i}$.
Let $\Gamma= C_2 \oplus C_2 = \langle \gamma\rangle \oplus \langle \eta\rangle$
with the action of $A$ given by $a(\gamma)=\gamma$ and $a(\eta)=\gamma+\eta$.
Then $\F\Gamma$ is a color Hopf algebra since $\beta(x,y)=1$ for all $x,y\in\supp \F\Gamma =\{1,g^2\}$.
\end{exa}

\begin{exa}\label{colorgroup2} Let $G= C_2\oplus C_4 = \langle g\rangle\oplus \langle h\rangle$,
${\mathbf i}\in \F$ such that ${\mathbf i}^2=-1$ and $\beta:G\times G \to \F^{\times}$, $\beta(g^i+h^j,g^k+h^l)={\mathbf i}^{il-jk}$;
then $A=\langle a\rangle\oplus \langle b\rangle$, where $a(g)=-1$, $a(h)=b(g)=1$ and $b(h)={\mathbf i}$.
Let $\Gamma= C_2 \oplus C_2 = \langle \gamma\rangle \oplus \langle \eta\rangle$
with the action of $A$ given by $a(\gamma)=b(\gamma)=\gamma$ and $a(\eta)=b(\eta)=\gamma+\eta$.
Then $\F\Gamma$ is a color Hopf algebra since $\beta(x,y)=1$ for all $x,y\in\supp \F\Gamma =\{1,g+h^2\}$.
\end{exa}

\subsection{Abelian extensions in color categories}

In this subsection, we analyze color Hopf algebras that arise as abelian extensions.

\subsubsection{Abelian extensions}
Let $(\Gamma,L)$ be a matched pair of groups \cite{Ma,Ta}, that is, $L$ and $\Gamma$ are finite groups provided with (right and left) actions
$\lhd:L\times \Gamma\to L$ and $\rhd: L\times \Gamma\to \Gamma$ satisfying: for all $l,t\in L$ and $\gamma,\eta\in \Gamma$
$$l\rhd \gamma\eta=(l\rhd \gamma)((l\lhd \gamma)\rhd\eta)\,\,\,\,\text{and}\,\,\,\,lt\lhd \gamma=(l\lhd (t\rhd \gamma))(t\lhd \gamma).$$

\nd Let $\sigma:\Gamma\times \Gamma\to (\F^{\times})^L$ and $\tau:L\times L\to (\F^{\times})^{\Gamma}$ be normalized 2-cocycles;
write $\sigma=\sum_{l\in L}\sigma_l\delta_l$,  $\tau=\sum_{\gamma\in \Gamma}\tau_{\gamma}\delta_{\gamma}$.
We shall always assume that $\sigma$ and $\tau$ satisfy the following normalization conditions: for all $l,t\in L$ and $\gamma,\eta\in\Gamma$
\begin{equation}\label{normal}
\sigma_1(\gamma,\eta)=1,\,\,\,\tau_1(l,t)=1.
\end{equation}

Let $H = \F^{L} {\,}^{\tau}\hspace{-0.08cm}\#_{\sigma} \F\Gamma$ be the vector space $\F^{L}\otimes \F\Gamma$ with the crossed product algebra and
the crossed coproduct coalgebra structures. The multiplication and coproduct of $H$ are
\begin{align}\label{eqcoc1}
&(\delta_le_{\gamma})(\delta_{t}e_{\eta}) = \delta_{l\lhd\gamma,t}\sigma_l(\gamma,\eta)\delta_le_{\gamma\eta}, & &l,t\in L, & &\gamma,\eta \in \Gamma,
\\ \label{eqcoc2}
&\Delta(\delta_le_{\gamma}) = \sum\limits_{u\in L}\tau_{\gamma}(u,u^{-1}l)\delta_{u}e_{(u^{-1}l)\rhd\gamma}\otimes\delta_{u^{-1}l}e_{\gamma}, & &l \in L, &
&\gamma \in \Gamma.
\end{align}
Here $\delta_le_{\gamma} := \delta_l\otimes e_{\gamma}\in H$. Then the following holds:

\smallbreak
\nd $\bullet$ \cite{Kac,Si,Ta}
$H$ is a Hopf algebra  iff  for all $x,y\in \Gamma$, $s,t\in L$.
\begin{equation}\label{eq:co44}
\sigma_{st}(x,y)\tau_{xy}(s,t) = \sigma_{s}(t\trid x,(t\trii x)\trid y)
	\sigma_{t}(x,y)\tau_{x}(s,t)\tau_{y}(s\trii(t\trid x),t\trii x).
\end{equation}
This happens e.~g. if $\sigma$ and $\tau$ are trivial, case where we denote  $H = \F^{L}\# \F\Gamma$.
If \eqref{eq:co44} holds, then we have an exact sequence of Hopf algebras
\begin{equation}\label{eqex4}
\xymatrix{1\ar[r] & \F^{L} \ar^\iota@{^{(}->}[r] & H \ar^\pi@{->>}[r] & \F\Gamma \ar[r] & 1.}
\end{equation}

\smallbreak
\nd $\bullet$ \cite[Theorem 2.10]{AN} gives necessary and sufficient conditions on $\sigma$ and $\tau$ so that $H = \F^{L} {\,}^{\tau}\hspace{-0.08cm}\#_{\sigma} \F\Gamma$ becomes a braided Hopf algebra with respect to a braiding $c: H\otimes H \to H \otimes H$. This braiding $c$ is uniquely determined by $\sigma$ and $\tau$, and is diagonal with respect to the basis $(\delta_le_{\gamma})_{l \in L, \gamma \in \Gamma}$. If these conditions hold, then $H$ can be realized as a braided Hopf algebra over an abelian group. We discuss this further in Subsection \ref{subsec:ext-nontrivial-braiding}.

\subsubsection{Automorphisms of an abelian extension with trivial cocycles}\label{subsubsection:aut-ext-trivial-cocycles}
In subsection \ref{subsubsection:color-ext-trivial-cocycles}, we shall assume that $\sigma,\tau$ are trivial 2-cocycles. We shall
give conditions on $H = \F^{L} \# \F\Gamma$ to be a $(G, \beta)$-color Hopf algebra, under the assumption that
$A < \Aut_{\text{ext}} H :=\{f\in \Aut H\,:\,f(x)\in \F^L,\,\,\text{for all}\,\,x\in \F^L\}$.
In this subsection, we determine $\Aut_{\text{ext}} H$, $H = \F^{L} \# \F\Gamma$.
We first observe that this group is sometimes the full automorphism group.

\begin{lema} Then the following are equivalent:
\begin{enumerate}\renewcommand{\theenumi}{\roman{enumi}}   \renewcommand{\labelenumi}{(\theenumi)}
\item $\mathcal{HZ}(H)=\F^L$,
\item $\lhd$ is trivial and $Z(\Gamma)=\{1\}$.
\end{enumerate}
\nd If either of these conditions holds, then $\Aut H=\Aut_{\textrm{ext}} H$.
\end{lema}

\pf (i) $\Rightarrow$ (ii): Since $\delta_{l,t}\delta_le_{\gamma}=\delta_l(\delta_te_{\gamma})=(\delta_te_{\gamma})\delta_l=\delta_{t\lhd\gamma,l}\delta_te_{\gamma}$, for all $l,t\in L$ and $\gamma\in \Gamma$,  $\lhd$ is trivial. The subalgebra of $H$ generated by $\{\delta_le_{\gamma}\,:\,l\in L,\,\,\gamma\in Z(\Gamma)\}$ is a central Hopf subalgebra of $H$, hence $Z(\Gamma)=\{1\}$.

\nd (ii) $\Rightarrow$ (i): Clearly,  $\pi(\mathcal{HZ}(H)) = \F$, because the former is a central Hopf subalgebra of $\F\Gamma$. Hence  $\pi(y) = \varepsilon(y)1$ for any $y\in \mathcal{HZ}(H)$. Therefore, $\mathcal{HZ}(H)\subseteq LKer(\pi)=\F^L$. Since $\lhd$ is trivial, $\F^L\subseteq \mathcal{HZ}(H)$.
The last assertion follows because the Hopf center is invariant under $\Aut H$.
\epf

Let $f\in \Aut_{\text{ext}} H$. If $l\in L$, then $f(\delta_l)=\delta_{f_1(l)}$, for some $f_1(l)\in L$;
since $f$ is a automorphism of Hopf algebras, we see that $f$ induces $f_1\in \Aut L$.
Analogously, $f$ induces  $f_2\in \Aut\Gamma$, and we have a homomorphism of groups
\begin{align*}
\Phi:\Aut_{\text{ext}} H &\to \Aut L\times\Aut\Gamma, & \Phi(f) & =(f_1,f_2).
\end{align*}

\begin{prop}\label{prop:3.8} Let $(g,h)\in \Aut L\times \Aut\Gamma$. There exists $f\in \Aut_{\text{\rm ext}} H$
such that $\Phi(f)=(g,h)$ if and only if the following conditions are satisfied:
\begin{enumerate}\renewcommand{\theenumi}{\roman{enumi}}   \renewcommand{\labelenumi}{(\theenumi)}
\item $g(l)\lhd h(\gamma)=g(l\lhd \gamma)$ and $g(l)\rhd h(\gamma)=h(l\rhd \gamma),\,\,\,\,l\in L,\,\,\gamma\in \Gamma$.\vu

\item  There exists a  map $\widetilde{f}: \Gamma\to (\F^{\times})^L$,
$\gamma\mapsto \widetilde{f}_{\gamma}:L\to \F^{\times}$ (uniquely determined, see \eqref{eq:f-tilde}) such that

\begin{align}\label{cyc1}
\widetilde{f}_{\gamma}(1) &= 1, & &\gamma\in\Gamma,
\\\label{cyc2}
\widetilde{f}_1(l) &= 1, & &l\in L,
\\\label{cyc3}
\widetilde{f}_{\gamma\eta}(l) &= \widetilde{f}_{\gamma}(l)\widetilde{f}_{\eta}(l\lhd h(\gamma)), & &l\in L,& &\gamma,\eta \in \Gamma,
\\\label{cyc4}
\widetilde{f}_{\gamma}(lt) &= \widetilde{f}_{g^{-1}(t)\rhd \gamma}(l)\widetilde{f}_{\gamma}(t),& &l,t\in L,& &\gamma\in \Gamma.
\end{align}
\end{enumerate}
\nd If (i) and (ii) hold, then $f$ is given by
\begin{align}\label{eq:def-f}
f(\delta_le_{\gamma}) &= \widetilde{f}_{\gamma}(g(l))\delta_{g(l)}e_{h(\gamma)},  &l\in L, \quad\gamma\in \Gamma.
\end{align}

\end{prop}

\pf Let $f\in \Aut_{\text{ext}} H$ such that $\Phi(f)=(g,h)$. Set  $\widetilde{e}_{\gamma} = 1 \otimes e_{\gamma}$, $\gamma\in \Gamma$.
Note that $f(\delta_le_{\gamma})=
f(\delta_le_1\cdot \widetilde{e}_{\gamma})=\delta_{g(l)}f(\widetilde{e}_{\gamma})$, for all $l\in L$ and $\gamma\in \Gamma$. Write
\begin{equation}\label{equation1}f(\widetilde{e}_{\gamma})
=\sum\limits_{\eta\in \Gamma}\sum\limits_{t\in L} \widetilde{f}_{\gamma}(t,\eta)\delta_te_{\eta},\,\,\,\widetilde{f}_{\gamma}(t,\eta)\in \F.\end{equation}
Then $\sum\limits_{\eta\in \Gamma}\widetilde{f}_{\gamma}(1,\eta)e_{\eta}=\pi\circ f(\widetilde{e}_{\gamma})=h(e_{\gamma})=e_{h(\gamma)}$, and consequently,
\begin{equation}\label{equation2} \widetilde{f}_{\gamma}(1,\eta)=\delta_{h(\gamma),\eta}.\end{equation}
Given $\gamma,\tau\in \Gamma$, we have that

\begin{align*}
f(\widetilde{e}_{\gamma})f(\widetilde{e}_{\tau})&=
\left(\sum\limits_{u\in \Gamma}\sum\limits_{m\in L}\widetilde{f}_{\gamma}(m,u)\delta_me_{u}\right) \left(\sum\limits_{v\in \Gamma}\sum\limits_{n\in L}\widetilde{f}_{\tau}(n,v)\delta_ne_{v}\right)\\
&= \sum\limits_{u,v\in \Gamma}\sum\limits_{l\in L}\widetilde{f}_{\gamma}(n\lhd u^{-1},u)\widetilde{f}_{\tau}(n,v)\delta_{n\lhd u^{-1}}
e_{u v};
\\
f(\widetilde{e}_{\gamma\tau}) &=
\sum\limits_{\eta\in \Gamma}\sum\limits_{l\in L}\widetilde{f}_{\gamma\tau}(l,\eta)\delta_{l}e_{\eta}\end{align*}
Hence
\begin{equation}\label{equation3} \widetilde{f}_{\gamma\tau}(l,\eta)=\sum\limits_{v\in \Gamma}\widetilde{f}_{\gamma}(l,\eta v^{-1})\widetilde{f}_{\tau}(l\lhd\eta v^{-1},v).
\end{equation}
 Also, $\Delta\circ f(\widetilde{e}_{\gamma}) = \sum\limits_{\eta\in \Gamma}\sum\limits_{r,t\in L}\widetilde{f}_{\gamma}(t,\eta)
\delta_{r}e_{(r^{-1}t)\rhd \eta}\otimes\delta_{r^{-1}t}e_{\eta}$ and
\begin{align*}
(f\otimes f)\circ \Delta(\widetilde{e}_{\gamma})&= (f\otimes f)\left(\sum\limits_{l\in L}\Delta(\delta_le_{\gamma})\right)
=\sum\limits_{l,s\in L}f(\delta_se_{(s^{-1}l)\rhd\gamma})\otimes f(\delta_{s^{-1}l}e_{\gamma})\\
&=\sum\limits_{l,s\in L}\delta_{g(s)}f(1_{\F^L}e_{(s^{-1}l)\rhd\gamma})\otimes \delta_{g(s^{-1}l)}f(\widetilde{e}_{\gamma})
\\
\overset{\eqref{equation1}}= & \sum\limits_{\substack{u,v\in \Gamma\\l,s\in L}}\widetilde{f}_{(s^{-1}l)\rhd\gamma}(g(s),u)
\widetilde{f}_{\gamma}(g(s^{-1}l),v)   \delta_{g(s)}e_{u}\otimes\delta_{g(s^{-1}l)}e_{v}.
\end{align*}
This implies that
\begin{equation}\label{eqcom}
\widetilde{f}_{\gamma}(rt,\eta)=\widetilde{f}_{g^{-1}(t)\rhd\gamma}(r,t\rhd \eta)\widetilde{f}_{\gamma}(t,\eta).
\end{equation}
It follows from \eqref{equation2} and \eqref{eqcom} that $\widetilde{f}_{\gamma}(r,\eta)=0$ whenever $\eta\neq h(\gamma)$.
Consider the map $\widetilde{f}:\Gamma\to (\F^{\times})^L,\,\,\gamma\mapsto \widetilde{f}_{\gamma}$, given by
\begin{align}\label{eq:f-tilde}
\widetilde{f}_{\gamma}(l) & = \widetilde{f}_{\gamma}(l,h(\gamma)), & &l\in L.
\end{align}
Then \eqref{cyc1} and \eqref{cyc3} follow, respectively, by \eqref{equation2} and \eqref{equation3}. Also, \eqref{cyc2} follows since $\sum\limits_{l\in L}\delta_le_{1}=1_{\F^L}e_1=f(1_{\F^L}e_1)=\sum\limits_{l\in L}\widetilde{f}_{1}(l)\delta_le_{1}$.

\nd Using that $f(\delta_le_{\gamma})=\delta_{g(l)}f(\widetilde{e}_{\gamma})=\widetilde{f}_{\gamma}(g(l))\delta_{g(l)}e_{h(\gamma)}$, we obtain
\begin{align}\label{eqcom1}
f(\delta_le_{\gamma}\cdot \delta_{t}e_{\eta}) &= \delta_{l\lhd \gamma,t}\widetilde{f}_{\gamma\eta}(g(l))\delta_{g(l)}e_{h(\gamma\eta)},
\\\label{eqcom2}
f(\delta_le_{\gamma})f(\delta_te_{\eta}) &= \widetilde{f}_{\gamma}(g(l))\widetilde{f}_{\eta}(g(t))\delta_{g(l)\lhd h(\gamma),g(t)}\delta_{g(l)}e_{h(\gamma)h(\eta)}.
\end{align}
The identity $g(l)\lhd h(\gamma)=g(l\lhd \gamma)$ is immediate from \eqref{eqcom1} and \eqref{eqcom2}. Using that $(f\otimes f)\circ \Delta(\delta_le_{\gamma})=\Delta\circ f(\delta_le_{\gamma})$, we obtain that $g(l)\rhd h(\gamma)=h(l\rhd \gamma)$. Finally, the last identity and \eqref{eqcom} imply \eqref{cyc4}.

\nd Conversely, assume that (i) and (ii) hold. Define $f$ by \eqref{eq:def-f};
by a straightforward calculation, $f\in \Aut_{\text{ext}} H$ and $\Phi(f)=(g,h)$.
\epf

\begin{remark}
The conditions in Proposition \ref{prop:3.8} (ii) can be spelled out in cohomology terms.
Indeed,  the action $\lhd$ induces actions of $\Gamma$ on $(\F^{\times})^L$ parameterized by $h\in\Aut\Gamma$:
$(\gamma\rightharpoonup_h \phi)\, (l)=\phi(l\lhd h(\gamma))$, $\gamma\in \Gamma$, $\phi\in (\F^{\times})^L$ and $l\in L$. Then \eqref{cyc2} and \eqref{cyc3} imply that $\widetilde{f}$ is a normalized 1-cocycle. When $\rhd$ is trivial, \eqref{cyc1} and \eqref{cyc4} imply that $\widetilde{f}\in Z^1(\Gamma,\Hom(L,\F^{\times}))$.
\end{remark}

\begin{exa}\label{forpro:3.7} Let $\Gamma= C_3 =\langle\gamma\rangle$ and $L=C_7 =\langle l\rangle$.
Then $(\Gamma,L)$ is a matched pair with trivial $\rhd$ and $\lhd:L\times \Gamma\to L$,
$l\lhd \gamma=l^{2}$. Then $g\in \Aut L$, $g(l)=l^{-1}$, and $h= \id\in \Aut \Gamma$,  satisfy the condition
Proposition \ref{prop:3.8} (i). Let $1\neq\xi\in\F$ such that $\xi^7=1$. The map $\widetilde{f}:\Gamma\to \F^{L}$,
$\widetilde{f}_1=1$, $\widetilde{f}_{\gamma}(l)=\xi$ and $\widetilde{f}_{\gamma^2}(l)=\xi^{3}$,
satisfies \eqref{cyc1}, \eqref{cyc2}, \eqref{cyc3} and \eqref{cyc4}. By Proposition \ref{prop:3.8}, there exists $f\in \Aut_{\text{ext}} H$
such that $\Phi(f)=(g,h)$. By \eqref{eq:def-f}, $f$ is given by
\begin{align*}
f(\delta_{l^i}e_1)&=\delta_{l^{-i}}e_1,&
f(\delta_{l^i}e_{\gamma})&=\xi^{-i}\delta_{l^{-i}}e_{\gamma},&
f(\delta_{l^i}e_{\gamma^2})&=\xi^{-3i}\delta_{l^{-i}}e_{\gamma^2}.
\end{align*}
\end{exa}

\begin{exa}\label{forpro1:3.7} Let $\Gamma= C_3 =\langle\gamma\rangle$ and $L=C_{12} =\langle l\rangle$. Then $(\Gamma,L)$ is a matched pair with the following actions
\begin{align*}
&l^{2k+1}\rhd \gamma =l^{2k}\rhd \gamma^2=\gamma^2,&  &l^{2k+1}\lhd \gamma = l^{2k+5},& l^{2k+1}\lhd \gamma^2&=l^{2k+9},&\\
&l^{2k}\rhd \gamma =l^{2k+1}\rhd \gamma^2=\gamma,& &l^{2k}\lhd \gamma=l^{2k}\lhd \gamma^2=l^{2k},& l^j\rhd 1&=1.
\end{align*}
Also $g\in \Aut L$, $g(l)=l^{7}$, and $h= \id\in \Aut \Gamma$,  satisfy the condition
Proposition \ref{prop:3.8} (i). Let $1\neq\xi\in\F$ such that $\xi^3=1$. Then $\widetilde{f}:\Gamma\to \F^{L}$, $\widetilde{f}_1=1$, $\widetilde{f}_{\gamma}(l^{2k})=\widetilde{f}_{\gamma^2}(l^{2k})=1$, $\widetilde{f}_{\gamma}(l^{2k+1})=\xi$ and $\widetilde{f}_{\gamma^2}(l^{2k+1})=\xi^{2}$,
satisfies \eqref{cyc1}, \eqref{cyc2}, \eqref{cyc3} and \eqref{cyc4}. By Proposition \ref{prop:3.8}, there exists $f\in \Aut_{\text{ext}} H$
such that $\Phi(f)=(g,h)$. By \eqref{eq:def-f}, $f$ is given by
\begin{align*}
f(\delta_{l^{2k}}e_{\gamma^i})&=\delta_{l^{2k}}e_{\gamma^i},&
f(\delta_{l^{2k+1}}e_{1})&=\delta_{l^{2(k+3)+1}}e_{1},&\\
f(\delta_{l^{2k+1}}e_{\gamma})&=\xi\delta_{l^{2(k+3)+1}}e_{\gamma},&
f(\delta_{l^{2k+1}}e_{\gamma^2})&=\xi^2\delta_{l^{2(k+3)+1}}e_{\gamma^2}.
\end{align*}
\end{exa}

\subsubsection{Color abelian extension}\label{subsubsection:color-ext-trivial-cocycles}

Let $\rho:A \to \Aut_{\text{ext}} H$ be a morphism of groups; we investigate when $H$ is a $(G, \beta)$-color Hopf algebra.

Given $a\in A$, we recall that there exists $(a_1,a_2) \in \Aut L\times\Aut\Gamma$
and $\tilde{a}_{\gamma}: \Gamma\to (\F^{\times})^L$
such that $\Phi(a)=(a_1,a_2)$ and $a(\delta_le_{\gamma})=\tilde{a}_{\gamma}(a_1(l))\delta_{a_1(l)}e_{a_2(\gamma)}$, $l\in L$
and $\gamma\in \Gamma$; see Proposition \ref{prop:3.8}.  We denote
$$A^{l}_{\gamma}:=\{a\in A\,:a_1(l)=l,\,\,a_2(\gamma)=\gamma\},\,\,\,\,G^{l}_{\gamma}:=\{g\in G\,:\,\chi_g\in A^{l}_{\gamma}\}.$$

Clearly, the map $\varphi^{l}_{\gamma}\,:\,A^{l}_{\gamma}\to \F^{\times}$, $\varphi^{l}_{\gamma}(a)=\tilde{a}_{\gamma}(l)$,
is a character of $A^{l}_{\gamma}$. Let $\chi^l_{\gamma}\,:\,G^{l}_{\gamma}\to \F^{\times}$, $\chi^l_{\gamma}=\varphi^{l}_{\gamma}\circ\chi$,
and $\rho^l_{\gamma}:A\to \Hom(G^l_{\gamma},\F^{\times})$, $\rho^l_{\gamma}(a)=a|_{G^l_{\gamma}}$.

\begin{definition}\label{def:3.11} A {\it $(G,\beta)$-color matched pair} is a matched pair $(\Gamma,L)$ provided with
a morphism of groups $\rho:A \to \Aut_{\text{ext}} H$ that satisfies: for all $u,l\in L$ and $\eta,\gamma\in \Gamma$,
\end{definition}
\begin{enumerate}\renewcommand{\theenumi}{\roman{enumi}}   \renewcommand{\labelenumi}{(\theenumi)}
\item $(G^{u}_{\eta})^{\bot}< A^l_{\gamma}$,
\item $A^l_{\gamma}\cap (\rho^u_{\eta})^{-1}(\chi^u_{\eta})\neq\emptyset$,
\item if $a\in A^l_{\gamma}\cap (\rho^u_{\eta})^{-1}(\chi^u_{\eta})$, then $\tilde{a}_{\gamma}(l)=1$.
\end{enumerate}

\begin{theorem}\label{th:3.13} Let $(\Gamma,L)$ be a matched pair provided with
a morphism of groups $\rho:A \to \Aut_{\text{ext}} H$.  Then the following are equivalent:
\begin{enumerate}\renewcommand{\theenumi}{\roman{enumi}}   \renewcommand{\labelenumi}{(\theenumi)}
\item $H$ is a color Hopf algebra,
\item $\beta(g,h)=1$, for all $g,h\in \text{sup}\,H$,
\item  $(\Gamma,L)$ is a $(G,\beta)$-color matched pair.
\end{enumerate}
\end{theorem}

\pf (i) $\Leftrightarrow$ (ii) follows from Proposition \ref{prop:sch}.
 (ii) $\Rightarrow$ (iii): Let $l,u\in L$, $\gamma,\eta\in \Gamma$ and $a\in (G^u_{\eta})^{\bot}$. Since $\text{sup}\,H\subset G^{u}_{\eta}$, $a\in A^{l}_{\gamma}$.
Let $b\in A$ such that $b|_{G^u_{\eta}}=\chi^u_{\eta}$ and fix $x=\delta_le_{\gamma}=\sum_{g\in s(x)}x_g$. Since $s(x)\subset \text{sup}\,H\subset G^u_{\eta}$, it follows that $b(g)=(\tilde{\chi}_g)_{\eta}(u)=1$, for all $g\in s(x)$. Then, $x=b(x)$ and  consequently, $b\in A^l_{\gamma}\cap (\rho^u_{\eta})^{-1}(\chi^u_{\eta})$. Finally, if $c\in A^l_{\gamma}\cap (\rho^u_{\eta})^{-1}(\chi^u_{\eta})$, then $c(x)=x$. Hence, $\tilde{c}_{\gamma}(l)=1$.
 (iii) $\Rightarrow$ (ii) Let $g,h\in \text{sup}\,H$. Consider $l,u\in L$ and $\gamma,\eta\in \Gamma$ such that $g\in s(\delta_ue_{\eta})$ and $h\in s(\delta_le_{\gamma})$. If $k\in G^{u}_{\eta}$, then $\chi_{g^{-1}}(k)=\chi^u_{\eta}(k)$, that is, $\rho^{u}_{\eta}(\chi_{g^{-1}})=\chi^u_{\eta}$.
Let $a\in A^l_{\gamma}$ with $\rho^u_{\eta}(a)=\chi^u_{\eta}$. Then, $a\chi_g\in (G^u_{\eta})^{\bot}$, which implies that $\chi_{g^{-1}}\in A^l_{\gamma}$.
Thus, $1=(\tilde{\chi}_{g^{-1}})_{\gamma}(l)=\chi_{g^{-1}}(h)=\beta(g,h)^{-1}$.
\epf

\begin{exa}\label{matchedpair} Let $H$ and $f$ be as in Example \ref{forpro1:3.7}, $G= C_2\oplus C_2 = \langle g\rangle\oplus \langle h\rangle$ and $\beta:G\times G \to \F^{\times}$, $\beta(g^i+h^j,g^k+h^l)=(-1)^{il-jk}$.
Then $A=\langle a\rangle\oplus \langle b\rangle$, where $a(g)=-1$, $a(h)=b(g)=1$ and $b(h)=-1$.
Consider the action of $A$ on $H$ given by $a(u)=b(u)=f(u)$, $ab(u)=u$, for all $u\in H$. Then
$H$ is a color Hopf algebra since $\beta(x,y)=1$ for all $x,y\in\supp H =\{1,g+h\}$.
\end{exa}

\subsection{Color abelian extensions with non-trivial braiding}\label{subsec:ext-nontrivial-braiding}

\

\nd Let $\sigma:\Gamma\times \Gamma\to (\F^{\times})^L$ and $\tau:L\times L\to (\F^{\times})^{\Gamma}$ be normalized 2-cocycles satisfying \eqref{normal} and let $H=\F^{L} {\,}^{\tau}\hspace{-0.08cm}\#_{\sigma} \F\Gamma$ with product \eqref{eqcoc1} and coproduct \eqref{eqcoc2}. In \cite[Section 3]{AN}, the authors discuss when $\F^L {\,}^{\tau}\hspace{-0.08cm}\#_{\sigma} \F\Gamma$ is a Hopf algebra in $^{\F G}_{\F G}\mathcal{YD}$ with respect to an action and a coaction
that are diagonal with respect to the basis $(\delta_le_{\gamma})_{l \in L, \gamma \in \Gamma}$; namely, they are determined  by
 maps $z:L\times\Gamma\to G$ and $\omega:L\times\Gamma\to \Hom(G,\F^{\times})$,  by the rules
\begin{align}\label{actioncoation}
g\cdot\delta_le_{\gamma} &=\omega(l,\gamma)(g)\delta_le_{\gamma}, & \lambda(\delta_le_{\gamma}) &= z(l,\gamma)\otimes\delta_le_{\gamma}.
\end{align}

Since we are interested in the braided tensor category $(\Vect^G, c_\beta)$, it is enough to
fix $z$ because $\omega$
is given by $\omega(l,\gamma)(g)=\beta(g,z(l,\gamma))$, see \eqref{eq:explicit-tensorfunctor}.
Note that \eqref{eqcoc1} and \eqref{eqcoc2} are morphisms of $\F G$-comodules if and only if,
\begin{align}\label{eqle1}
z(l,\gamma\eta)&=z(l,\gamma)z(l\lhd \gamma,\eta),& l&\in L,& \gamma,\eta&\in\Gamma,
\\\label{eqle2}
z(lt,\gamma)&=z(l,t\rhd\gamma)z(t,\gamma),& l,t&\in L, & \gamma\in \Gamma.
\end{align}

The following theorem is a consequence of \cite[Theorem 3.5]{AN}. \vu

\begin{theorem}\label{thm:3.3.15}  Suppose that $z:L\times\Gamma\to G$ satisfies
\eqref{eqle1} and \eqref{eqle2}. Then $H$ is a $(G,\beta)$-color Hopf algebra iff for all $l,t\in L$ and $\gamma,\eta\in\Gamma$,
\begin{align}\label{eq:compatibility}
\begin{aligned}
&\sigma_{lt}(\gamma,\eta)\tau_{\gamma\eta}(l,t)  =\beta\left(z(t,\gamma), \, z\big(l\lhd(t\rhd \gamma),(t\lhd\gamma)\rhd \eta\big)\right)\,
\\ &\times \tau_{\gamma}(l,t) \,
\tau_{\eta}(l\lhd(t\rhd\gamma),t\lhd\gamma) \,
\sigma_l(t\rhd\gamma,(t\lhd\gamma)\rhd\eta) \,\sigma_t(\gamma,\eta). \qed
\end{aligned}
\end{align}
\end{theorem}

For the rest of this subsection, in order to determine maps $z$ satisfying \eqref{eqle1} and \eqref{eqle2}, we assume that $\rhd$ is trivial; this implies that $\lhd:L\times \Gamma\to L$ is an action by group automorphisms.
Then $\lhd$ induces left actions of $\Gamma$ on $\Hom(L,G)$ and on $\Hom(L,A)$; for instance, $(\gamma\rightharpoonup \phi)(l)=\phi(l\lhd \gamma)$, for all $\gamma\in \Gamma$, $\phi\in \Hom(L,G)$ and $l\in L$. By \cite[Lemma 4.1]{AN}, the correspondence
$z\mapsto \widetilde{z}:\Gamma\to \Hom(L,G)$, $\widetilde{z}(\gamma)(l)=z(l,\gamma)$, $l\in L$, $\gamma\in \Gamma$,
determines a bijection between the set of maps $z$ satisfying \eqref{eqle1} and \eqref{eqle2} and $Z^1(\Gamma,\Hom(L,G))$.
We now look for $\sigma,\tau$ satisfying the compatibility condition \eqref{eq:compatibility}. We consider the 2-cocycle $\tau$ as a map $\widetilde{\tau}:\Gamma\to Z^2(L,\F^{\times})$ and the  action of $\Gamma$ on $Z^2(L,\F^{\times})$:
$(\gamma\cdot\varphi)(l,t)=\varphi(l\lhd \gamma,t\lhd\gamma)$, $\gamma\in \Gamma$, $\varphi\in Z^2(L,\F^{\times})$, $l,t\in L$.

\begin{theorem}\label{thm:3.17} Let $\widetilde{z}\in Z^1(\Gamma,\Hom(L,G))$  and suppose that the following compatibility condition holds
\begin{align}\label{eq:theorem3.17}
&\sigma_{lt}(\gamma,\eta)=\beta\left(\widetilde{z}(\gamma)(t),\gamma\rightharpoonup \widetilde{z}(\eta)(l)\right)\,\sigma_l(\gamma,\eta)\,\sigma_t(\gamma,\eta).&
\end{align}
Then $H$ is a $(G,\beta)$-color Hopf algebra if and only if $\widetilde{\tau}\in Z^1\left(\Gamma,Z^2(L,\F^{\times})\right)$.
\end{theorem}

\pf Immediate consequence of Theorem \ref{thm:3.3.15}.
\epf

We end this subsection with examples of color Hopf algebras that are neither commutative nor cocommutative. As in \cite[Chapter 3]{Sommer}
we consider:
\begin{align*}
R&  \text{ a finite ring},& L=G &= \text{ the additive group of } R,
\\
\Gamma&\text{ a finite group},& \nu:\Gamma\to& R^{\times} \text{ a group homomorphism,}
\\
\psi&\in Z^1(\Gamma,L), & \phi&\in  Z^2(\Gamma,L), \\
\eta&,\theta:L\to \F^{\times},& \theta(l&tu)=\theta(tlu), \forall l,t, u\in L,\\
l&\rhd \gamma=\gamma, \quad l\lhd\gamma=l\,\nu(\gamma),& l&\in L,\quad\gamma\in\Gamma.
\\
z&:L\times \Gamma\to G,\quad z(l,\gamma)=l\psi(\gamma),& l&\in L, \quad\gamma\in\Gamma, \\
\beta&:G\times G\to \F^{\times},\quad \beta(g,h)=\theta(gh)^2,&  g,h&\in G.
\end{align*}
Define $\sigma_l(\gamma,v)=\eta\left(l\phi(\gamma,v)\right)\theta\left(l^2\nu(\gamma)\psi(\gamma)\psi(v)\right)$, for
$l\in L$, $\gamma,v\in \Gamma$.
By \cite[3.3]{Sommer}, $\widetilde{z}$, $\beta$ and $\sigma$ satisfy \eqref{eq:theorem3.17}. By Theorem \ref{thm:3.17}, $H$ is a $(G,\beta)$-color Hopf algebra for each $\widetilde{\tau}\in Z^1(\Gamma,Z^2(L,\F^{\times}))$; see \cite[Lemma 4.7]{AN}
for examples of such 1-cocycles $\widetilde{\tau}$ of $\Gamma$ in $Z^2(L,\F^{\times})$.

\subsection*{Acknowledgments}
This paper was done during a posdoctoral stage of the third author in the
Universidad Nacional de C\'ordoba. He thanks N. Andruskiewitsch and I. Angiono for the kind hospitality.

\end{document}